\documentclass[12pt]{article}
\usepackage{bbm}
\usepackage{mathrsfs}
\usepackage{amsfonts}
\usepackage{amsmath,amssymb,latexsym,amsfonts}

\usepackage{CJK,color}
\usepackage{amssymb}
\usepackage{amsmath}
\usepackage{graphics}
\usepackage{epsfig}
\usepackage{indentfirst}
\usepackage{cite}

\definecolor{red}{rgb}{1.00,0.00,0.00}
\definecolor{blue}{rgb}{0.00,0.00,0.63}
\definecolor{black}{rgb}{0.00,0.00,0.00}

\newtheorem{claim}{\bf \t}[part]

\newtheorem{Lemma}{Lemma}[part]
\newtheorem{Proposition}{Proposition}[part]
\newtheorem{Remark}{Remark}[part]

\newtheorem{Theorem}{Theorem}[part]

\numberwithin{Assumption}{section} \numberwithin{Corollary}{section}
\numberwithin{Definition}{section} \numberwithin{equation}{section}
\numberwithin{Example}{section} \numberwithin{Lemma}{section}
\numberwithin{Proposition}{section} \numberwithin{Remark}{section}
\numberwithin{Theorem}{section}

\def\t{\theta}

\def\f{\frac}

\def \R{{\mathbb{R}}}

\def \qed{{\hbox{ }\hfill$\Box$}}

\def\text#1{{\rm #1}}

\begin{document}
%\begin{CJK*}{GBK}{song}
\date{}

\title{\Large \bf  Non-Existence of Classical Solutions with Finite Energy  to the Cauchy Problem of the Compressible Navier-Stokes  Equations }

\author{HAILIANG LI, YUEXUN WANG, AND ZHOUPING XIN}

%\small
\maketitle
\begin{abstract}
The well-posedness of classical solutions with finite energy to the compressible Navier-Stokes equations (CNS) subject to arbitrarily large and smooth initial data is a challenging problem. In the case when the fluid density is away from vacuum (strictly positive), this problem was first solved for the CNS in either one-dimension for general smooth initial data or multi-dimension for smooth initial data near some equilibrium state (i.e., small perturbation)~\cite{A,K1,K2,M1,M2,M3}. In the case that the flow density may contain vacuum (the density can be zero at some space-time point), it seems to be a rather subtle problem to deal with the well-posedness problem for CNS. The local well-posedness of classical solutions containing vacuum was shown in homogeneous Sobolev space (without the information of velocity in $L^2$-norm) for general regular initial data with some compatibility conditions being satisfied initially~\cite{CK1,CCK,CK2,CK3}, and the global existence of classical solution in the same space is established under additional assumption of small total initial energy but possible large oscillations~\cite{HLX}. However, it was shown that any classical solutions to the compressible Navier-Stokes equations in finite energy (inhomogeneous Sobolev) space can not exist globally in time since it may blow up in finite time provided that the density was compactly supported~\cite{Xin}. In this paper, we investigate the well-posedess of classical solutions to the Cauchy problem of Navier-Stokes equations,
 % in inhomogeneous Sobolev space. We
and prove that the classical solution with finite energy does not exist even in the inhomogeneous Sobolev space for any short time under some natural assumptions on initial data near the vacuum.  This implies in particular that the homogeneous Sobolev space is crucial as studying the well-posedness for the Cauchy problem of compressible  Navier-Stokes equations in the presence of vacuum at far fields even locally in time.
\end{abstract}

\setcounter{equation}{0}
\setcounter{Assumption}{0}
\setcounter{Theorem}{0}
\setcounter{Proposition}{0}
\setcounter{Corollary}{0}
\setcounter{Lemma}{0}
\mbox{}

\begin{center}
\section{Introduction and Main Results}
\end{center}

The motion of a $n$-dimensional compressible viscous, heat-conductive, Newtonian polytropic fluid
 is governed by the following full compressible
Navier-Stokes system:
\begin{eqnarray}\label{5}
\left\{ \begin{array}{ll}
\partial_t\rho+\textrm{div}(\rho
u)=0,\\
\partial_t(\rho u)+\textrm{div}(\rho u\otimes u)+\nabla p=\mu\Delta u+(\mu+\lambda)\nabla\textrm{div}u,\\
\partial_t(\rho e)+\textrm{div}(\rho eu)+ p\textrm{div}u=\frac{\mu}{2}|\nabla u+(\nabla u)^*|^2+\lambda(\textrm{div}u)^2+\frac{\kappa(\gamma-1) }{R}\Delta e,
\end{array}
\right.
\end{eqnarray}
where $(x,t)\in\mathbb{R}^n\times\mathbb{R}_+ $, $\rho,u,p$ and $e$ denote the
density, velocity, pressure and internal energy, respectively.  $\mu$ and $\lambda$ are the coefficient of viscosity
and the second coefficient of viscosity respectively and $\kappa$ denotes the coefficient of heat conduction, which satisfy
\begin{equation*}
 \mu>0,\quad  2\mu+n\lambda\geq0, \quad \kappa\geq0.
\end{equation*}
The  equation of state for polytropic gases satisfies
\begin{equation}\label{6}
 p=(\gamma-1)\rho e, \quad p=A\exp(\frac{(\gamma-1)S}{R})\rho^\gamma,
\end{equation}
where  $A>0$ and $R>0$ are positive constants, $\gamma> 1$  is the specific heat ratio, $S$ is the entropy, and we set $A=1$ in this paper for simplicity.   The initial data is given by
\begin{eqnarray}\label{7}
(\rho,u,e)(x,0)=(\rho_0,u_0,e_0)(x),\quad x\in\R^n
\end{eqnarray}
and is assumed to be continuous. In particular,  the initial density is compactly supported on an open bounded set  $\Omega\subset \mathbb{R}^n$  with smooth boundary, i.e.,
\begin{eqnarray}\label{8}
 \mbox{supp}_x\,\rho_0=\bar{\Omega},\quad  \rho_0(x)>0, \ x\in \Omega
\end{eqnarray}
and the initial internal energy $e_0$ is assumed to be nonnegative  but not identical to  zero in $\Omega$ to avoid the trivial case.

When the heat conduction can be neglected and the compressible viscous fluids are isentropic, the compressible Navier-Stokes equations \eqref{5} can be reduced to the following system
\begin{eqnarray}\label{1}
\left\{ \begin{array}{ll}
\partial_t\rho+\textrm{div}(\rho
u)=0,\\
\partial_t(\rho u)+\textrm{div}(\rho u\otimes u)+\nabla p=\mu\Delta u+(\mu+\lambda)\nabla\textrm{div}u,
\end{array}
\right.
\end{eqnarray}
for $(x,t)\in\mathbb{R}^n\times\mathbb{R}_+ $, where the  equation of state satisfies
\begin{equation}\label{2}
 p=A\rho^\gamma
\end{equation}
and the  initial data are given by
\begin{eqnarray}\label{3}
(\rho,u)(x,0)=(\rho_0,u_0)(x),\quad x\in\R^n
\end{eqnarray}
with the initial density being compactly supported, i.e., the assumption \eqref{8} holds.

It is an important issue to study the global existence (well-posedness) of classical/strong solution to CNS \eqref{5} and \eqref{1}, and many significant progress have been made recently on this and related topics, such as the global existence and asymptotical behaviors of solutions to \eqref{5} and \eqref{1}. For instance, in the case when the flow density is strictly away from the vacuum ($\inf_\Omega\rho>0$), the short time existence of classical solution was shown for general regular initial data~\cite{Ka}, the global existence of solutions problems were proved in spatial one-dimension by Kazhikhov et al. \cite{A,K1,K2} for sufficiently smooth data and by Serre \cite{S1,S2} and Hoff \cite{H1} for discontinuous initial data.
 The key point here behind the strategies to establish the global existence of strong solutions lies in the fact that if the flow density is strictly positive at the initial time, so does for any later-on time~\cite{H}.  This is also proved to be true for weak solutions to the compressible
Navier-Stokes equations (1.1) in one space dimension, namely, weak solution  does not exhibit vacuum states in any
finite time provided that no vacuum is present initially~\cite{H4}.
The corresponding multidimensional problems were also investigated as the flow density is away from the vacuum, for instance, the short time well-posedness of classical solution was shown by Nash and Serrin for general smooth initial data~\cite{Nash,Se}, and the global existence of unique strong solution was first proved by Matsumura and Nishida~\cite{M1,M2,M3} in the energy space (inhomogeneous Sobolev space)
\begin{equation}
\left\{\begin{aligned}\label{inhomSpace}
\rho-\bar{\rho}\in C(0,T;H^3(\mathbb{R}^3))\cap C^1(0,T;H^2(\mathbb{R}^3)), \\
u,\, e-\bar{e}\in C(0,T;H^3(\mathbb{R}^3))\cap C^1(0,T;H^1(\mathbb{R}^3)),
\end{aligned}\right.
\end{equation}
with $\bar{\rho}>0$ and $\bar{e}>0$ for any $T\in(0,\infty]$, where the additional assumption of small oscillation is required on the perturbation of initial data near the non-vacuum equilibrium state $(\bar{\rho},0,\bar{e})$. The global existence of non-vacuum solution was also solved by Hoff for discontinuous initial data \cite{H2}, and by Danchin \cite{D} who set up the framework based on the Besov type space (a functional space invariant by the natural scaling of the associated equations) to obtain existence and uniqueness of global solutions, where the small oscillations on the perturbation of initial data near some non-vacuum equilibrium state is also required.  It should be mentioned here that above smallness of the initial oscillation on the perturbation of initial data near the non-vacuum equilibrium state and the uniformly a-priori estimates established on the classical solutions  to CNS \eqref{5} or  \eqref{1} are sufficient to  establish the strict positivity and uniform bounds of flow density, which is essential to prove the global existence of solutions with the flow density away from vacuum in the inhomogeneous Sobolev space~\eqref{inhomSpace} or other function spaces~\cite{D,H2}. However, recently, this assumption on the small oscillations on the initial perturbation of a non-vacuum state can be removed at least for the isentropic case by Huang-Li-Xin in~\cite{HLX} provided that the initial total mechanical energy is suitable small which is equivalent to that the mean square norm of the initial difference from the non-vacuum state is small so that the perturbation may contain large oscillations and  vacuum state. See also~\cite{WenZhu}.

In the case when the flow density may contain vacuum (the flow density is nonnegative), it is rather difficult and challenging to investigate the global existence (well-posedness) of classical/strong solutions to CNS~\eqref{5} and CNS~\eqref{1}, corresponding  to the well-posedness theory of classical solutions~\cite{M1,M2,M3}, and the possible appearance of vacuum in the flow density (i.e., the flow density is zero)  is one of the essential difficulties in the analysis of the well-posedness and related problems~\cite{CCK,CK1,CK2,CK3,H1,H3,H4,Sa,S1,WenZhu,Xin,XYa,XYu}. Indeed, as it is well-known that \eqref{5} and \eqref{1} are strongly coupled systems of hyperbolic-parabolic type, the density $\rho(x,t)$ can be determined by its initial value $\rho_0(x_0)$ by Eq.~$\eqref{1}_1$ along the particle path $x(t)$ satisfying $x=x(t)$ and $x(0)=x_0$ provided that the flow velocity $u(x,t)$ is a-priorily regular enough.
Yet, the flow velocity can only be solved by Eq.~$\eqref{1}_2$ which is uniformly parabolic so long as the density is a-priorily strictly positive and uniformly bounded function. However, the appearance of vacuum leads to the strong degeneracy of the hyperbolic-parabolic system and the behaviors of the solution may become singular, such as the ill-posedness and finite blow-up of classical solutions~\cite{CJ,H3,S1,Xin,XYa}.
Recently, the global existence of weak solutions with finite energy to the isentropic system~\eqref{1} subject to general initial data with finite initial energy (initial data may include vacuum states) by Lions \cite{L1,L2,L3}, Jiang-Zhang~\cite{JZ} and Feireisl et al. \cite{F}, where  the exponent $\gamma$ may be required to be large and the flow density is allowed to vanish. Despite the important progress, the regularity, uniqueness and behavior of these
weak solutions remain largely open. As emphasized before~\cite{CJ,H3,S1,Xin,XYa}, the possible appearance of
vacuum is one of the major difficulties when trying to prove global existence and strong
regularity results. Indeed, Xin~\cite{Xin} first shows that it is impossible to obtain the global existence of finite energy classical solution to the Cauchy problem for \eqref{5} in the inhomogeneous Sobolev space~\eqref{inhomSpace} for any smooth initial data  with initial flow density compactly supported and similar phenomena happens for the isentropic system~\eqref{1} for a large class of smooth initial data with compactly supported density. To be more precise, if there exists any solution $(\rho,u,e)\in C^1(0,T;H^2(\mathbb{R}^3))$ for some time $T>0$, then it must hold $T<+\infty$, which also implies the finite time blow-up of solution $(\rho,u,e)\in C^1(0,T;H^2(\mathbb{R}^3))$ if existing in the presence of the vacuum.
Yet, Cho et al. \cite{CK1,CCK,CK2,CK3} proved the local well-posedness of classical solutions to the Cauchy problem for isentropic compressible Navier-Stokes equations~\eqref{1} and full Navier-Stokes equations~\eqref{5} with the initial density containing vacuum for some $T>0$ in the homogeneous energy space
\begin{equation}
\left\{\begin{aligned}\label{homSpace}
\rho\in C(0,T;H^3(\mathbb{R}^3))\cap C^1(0,T;H^2(\mathbb{R}^3)), \\
u,\, e\in C(0,T;D^3(\mathbb{R}^3))\cap L^2(0,T;D^{4}(\mathbb{R}^3)),
\end{aligned}\right.
\end{equation}
where $D^k(\mathbb{R}^3)=\{f\in L_{\textrm{{loc}}}^1(\mathbb{R}^3): \nabla f\in H^{k-1}(\mathbb{R}^3)\,\}$, under some additional compatibility conditions as \eqref{12.6}  on $u$ and similar compatibility condition on $e$.
Moreover, under additional smallness assumption on initial energy, the global existence and
uniqueness of classical solutions to the isentropic system~\eqref{1} established by Huang-Li-Xin in homogeneous Sobolev space~\cite{HLX}.  Interestingly, such a theory of global in time existence of classical solutions to the full CNS~\eqref{5} fails to be true due to the blow-up results Xin-Yan \cite{XYa} where they show that any classical solutions to ~\eqref{5} will blow-up in finite time as long as the initial density has an isolated mass group.  Note that the blow-up results in \cite{XYa} is independent of the spaces the solutions may be and whether they have small or large data.
It should be noted that the main difference of the homogeneous Sobolev space~\eqref{homSpace} from the inhomogeneous Sobolev space~\eqref{inhomSpace} lies that there is no any estimates on the term $\|u\|_{L^2}$ for the velocity.
Thus, it is natural and important to show whether or not the classical solution to the Cauchy problem for the CNS~\eqref{5} and CNS~\eqref{1} exits in the inhomogeneous Sobolev space~\eqref{inhomSpace} for some small time.

We study the well-posedess of classical solutions to the Cauchy problem for the full compressible Navier-Stokes equations~\eqref{5} and the isentropic Navier-Stokes equations~\eqref{1} in the  inhomogeneous Sobolev space~\eqref{inhomSpace} in the present paper, and we prove that there does not exist any classical solution in the  inhomogeneous Sobolev space~\eqref{inhomSpace} for any small time (refer to Theorems~\ref{1001}--\ref{1005} for details).  These imply that the homogeneous Sobolev spaces such as~\eqref{inhomSpace}, are crucial in the study of the well-posedness theory of classical solutions to the Cauchy problem of compressible Navier-Stokes equations in the presence of vacuum at far fields.

The main results in this paper can be stated as follows:
\begin{Theorem}
\label{1001}
The one-dimensional  isentropic Navier-Stokes  equations \eqref{1}-\eqref{3} with  the initial density satisfying \eqref{8} with $\Omega\triangleq I=(0,1)$ has no any solution $(\rho,u)$ in the inhomogeneous Sobolev space $C^1([0,T];H^m(\mathbb{R})), m>2$ for any positive time $T$, if the initial data $(\rho_0,u_0)$ satisfy one of the following two conditions in the interval $I$:
there exist  positive numbers $\lambda_i,i=1,2,3,4$ with $0<\lambda_3,\lambda_4<1$ such that
\begin{eqnarray}\label{11}
\left\{ \begin{array}{ll}
\frac{(\rho_0)_x}{\rho_0}\geq \lambda_1,\ in \ (0,\lambda_3),\\
u_0(\lambda_3)<0,u_0\leq0,\ in \ (0,\lambda_3),
\end{array}
\right.
\end{eqnarray}
or
\begin{eqnarray}\label{12}
\left\{ \begin{array}{ll}
\frac{(\rho_0)_x}{\rho_0}\leq-\lambda_2,\ in \ (\lambda_4,1),\\
u_0(\lambda_4)>0,u_0\geq0,\ in \ (\lambda_4,1).
\end{array}
\right.
\end{eqnarray}
\end{Theorem}

The following remark is helpful for understanding the conditions \eqref{11}-\eqref{12} and Theorem \ref{1001}.
\begin{Remark} The set of initial data $(\rho_0,u_0)$ satisfying the condition \eqref{11} or \eqref{12} is non-empty.
For example, for any given positive integers $k$ and $l$. Set
\begin{eqnarray}\label{12.2}
\rho_0(x)=
\left\{ \begin{array}{ll}
x^k(1-x)^k, \ for \ x\in \ [0,1],\\
0,\ for\ x\in \mathbb{R}\setminus [0,1]
\end{array}
\right.
\end{eqnarray}
and
\begin{eqnarray}\label{12.4}
u_0(x)=
\left\{ \begin{array}{ll}
- x^l, \ for\ x\in \ [0,\frac{1}{4}],\\
\ smooth\ connection, \ for\ x\in \ (\frac{1}{4},\frac{3}{4}),\\
 (1-x)^l, \ for\ x\in  \ [\frac{3}{4},1],\\
0, \ for\ x\in \mathbb{R}\setminus [0,1],
\end{array}
\right.
\end{eqnarray}
 then $(\rho_0,u_0)$ satisfies both \eqref{11} and \eqref{12}.

 It is known that the system \eqref{1}-\eqref{3} is well-poseded in the homogeneous Sobolev space in classical sense if and only if $\rho_0$ and $u_0$ satisfy the following compatibility condition (see \cite{CK2})
 \begin{eqnarray}\label{12.6}
\left\{ \begin{array}{ll}
-\mu\Delta u_0-(\mu+\lambda)\nabla\emph{div}u_0+\nabla p_0=\rho_0g,\\
g\in D^1,\sqrt{\rho_0}g\in L^2.
\end{array}
\right.
\end{eqnarray}
In one-dimensional case, for $(\rho_0,u_0)$ given by \eqref{12.2} and \eqref{12.4}, we have
\begin{eqnarray*}
g=
\left\{ \begin{array}{ll}
O(x^{l-k-2})+O(x^{l-k-1})+O(x^{k(\gamma-1)-1}), \ for \ x\in \ [0,\frac{1}{4}],\\
\ smooth\ connection, \ for\ x\in \ (\frac{1}{4},\frac{3}{4}),\\
O((1-x)^{l-k-2})+O((1-x)^{l-k-1})+O((1-x)^{k(\gamma-1)-1}),\ for\ x\in \ [\frac{3}{4},1],\\
0, \ for\ x\in \mathbb{R}\setminus [0,1].
\end{array}
\right.
\end{eqnarray*}
Direct calculations show $(\rho_0,u_0)$ satisfy \eqref{12.6} if and only if
\begin{eqnarray}\label{12.8}
\left\{ \begin{array}{ll}
k>\frac{3}{2(\gamma-1)},\\
l>k+\frac{5}{2}.
\end{array}
\right.
\end{eqnarray}

%Similarly, if we want to search solutions for the system \eqref{1}-\eqref{3} in  $ C^1([0,T];H^m(\mathbb{R}^n))$, $m>[\frac{n}{2}]+2$, then the
%compatibility condition should be
%\begin{eqnarray*}
%\left\{ \begin{array}{ll}
%-\mu\Delta u_0-(\mu+\lambda)\nabla\emph{div}u_0+\nabla p_0=\rho_0g,\\
%g\in D^1,\sqrt{\rho_0}g\in L^2,u_0\in L^2.
%\end{array}
%\right.
%\end{eqnarray*}

For the initial data $(\rho_0,u_0)$ given by  \eqref{12.2} and \eqref{12.4} with \eqref{12.8}, the system \eqref{1}-\eqref{3} is well-poseded in homogeneous Sobolev space but has no solution in  $ C^1([0,T];H^m(\mathbb{R}))$, $m>2$ for any  positive time $T$. Therefore, the solution constructed in (see \cite{CK2}) has no finite energy in $ C^1([0,T];H^m(\mathbb{R}))$, $m>2$  for any  positive time $T$ even if the initial data has finite energy in $H^m(\mathbb{R})$. Precisely, even if
\begin{eqnarray}
\int_{\mathbb{R}}u_0(x)^2dx<\infty,
\end{eqnarray}
but it holds that
\begin{eqnarray}
\int_{\mathbb{R}}u(x,t)^2dx=\infty \quad  \ for\ any\  t>0.
\end{eqnarray}
\end{Remark}

\begin{Theorem}\label{1003} The one-dimensional  full Navier-Stokes  equations  \eqref{5}-\eqref{7} with zero heat conduction and  the initial density satisfying \eqref{8} with $\Omega\triangleq I=(0,1)$ has no any solution $(\rho,u,e)$ in the inhomogeneous Sobolev space $C^1([0,T];H^m(\mathbb{R}))$, $m>2$  for any  positive time $T$, if the initial data $(\rho_0,u_0,e_0)$ satisfy one of the following two conditions in the interval $I$:
there exist  positive numbers $\lambda_i,i=5,6,7,8$ with $0<\lambda_7,\lambda_8<1$ such that
\begin{eqnarray}\label{14}
\left\{ \begin{array}{ll}
\frac{(\rho_0)_x}{\rho_0}+\frac{(e_0)_x}{\rho_0}\geq \lambda_5,\ in \ (0,\lambda_7),\\
u_0(\lambda_7)<0,u_0\leq0,\ in \ (0,\lambda_7),
\end{array}
\right.
\end{eqnarray}
or
\begin{eqnarray}\label{15}
\left\{ \begin{array}{ll}
\frac{(\rho_0)_x}{\rho_0}+\frac{(e_0)_x}{\rho_0}\leq-\lambda_6,\ in \ (\lambda_8,1),\\
u_0(\lambda_8)>0,u_0\geq0,\ in \ (\lambda_8,1).
\end{array}
\right.
\end{eqnarray}
\end{Theorem}

Huang and Li \cite{HL} proved the well-posedness  to the Cauchy problem of the $n$-dimensional full compressible Navier-Stokes equations \eqref{5}-\eqref{6} with positive heat conduction in Sobolev space, but the entropy function $S(t,x)$ is infinite in  vacuum domain (see Remark 4.2 in \cite{XYa}). If the entropy function $S(t,x)$ is required to be finite in  vacuum domains, then we have the following non-existence result:
\begin{Theorem}\label{1005}
The $n$-dimensional full compressible Navier-Stokes equations \eqref{5}-\eqref{7} with positive heat conduction
and  the initial density satisfying \eqref{8} has no any solution $(\rho,u,e)$ in the inhomogeneous Sobolev space $C^1([0,T];H^m(\mathbb{R}^n))$, $m>[\frac{n}{2}]+2$ with finite entropy $S(t,x)$ for any  positive time $T$.
\end{Theorem}

To prove  Theorem \ref{1001}-Theorem \ref{1005}, we will carry out the following steps.  First we reduce the original Cauchy problem to an initial-boundary value problem, which then can be reduced further to an integro-differential system  with degeneracy  for t-derivative  by the Lagrangian coordinates transformation, and one can then define a linear parabolic operator from the integro-differential system  and establish the Hopf's lemma and a strong maximum principle for the resulting  operator, and finally we prove that the resulting system is over-determined by contradiction.
Because the linear parabolic operator here degenerates for t-derivative due to that the initial density vanishes on boundary,  one needs careful analysis to deduce a localized version strong maximum principle on some rectangle away from boundaries.

We should stress that our method is based on maximum principle for parabolic operator, therefore we shall deal with one-dimensional isentropic case in Section 2, one-dimensional zero heat conduction case in Section 3 and n-dimensional positive heat conduction case in Section 4 separately, we define parabolic operators from momentum equation near the degenerate boundary in the Lagrangian coordinates  by adding some conditions on initial data for the first two cases and the energy equation in the whole domain for the last case, respectively.

\begin{center}
\section{Proof of Theorem \ref{1001}}
\end{center}

\subsection{Reformulation of Theorem \ref{1001}}

Let $n=1$ and $(\rho,u)\in C^1([0,T];H^m(\mathbb{R})), m>2$  be a solution  to the system \eqref{1}-\eqref{3} with  the initial density satisfying $\eqref{8}$.
Let $a(t)$ and $b(t)$ be the particle paths stating from $0$ and  $1$, respectively. The following argument is due to Xin \cite{Xin}.
Following  from the first equation of \eqref{1}, we see $\mbox{supp}_x\,\rho=[a(t), b(t)]$.
It follows from the second equation of \eqref{1} that
\begin{eqnarray*}
u_{xx}(x,t)=0, \forall \ x\in \mathbb{R}\backslash [a(t), b(t)]\times (0,T],
\end{eqnarray*}
which gives
\begin{eqnarray*}
u(x,t)=
\left\{ \begin{array}{ll}
u(b(t),t)+(x-b(t))u_x(b(t),t), \ if \ x>b(t),\\
u(a(t),t)+(x-a(t))u_x(a(t),t), \ if \ x<a(t).
\end{array}
\right.
\end{eqnarray*}
Since $u(\cdot,t)\in H^m(\mathbb{R}), m>2$,  then one has
\begin{eqnarray}\label{9}
u(x,t)=u_x(x,t)=0, \forall \ x\in \mathbb{R}\backslash [a(t), b(t)]\times (0,T],
\end{eqnarray}
which implies $[a(t), b(t)]=[0,1]$, i.e., $\mbox{supp}_x\,\rho(x,t)=[0,1]$.

Therefore, by the above argument,  to study the well-posedness of the system \eqref{1}-\eqref{3}  with  the initial density satisfying \eqref{8} is  equivalent to
study the well-posedness of the following  initial-boundary value problem
\begin{eqnarray}\label{10}
\left\{ \begin{array}{ll}
\rho_t+(\rho
u)_x=0,\ in\ I\times (0,T],\\
(\rho u)_t+(\rho u^2+p)_x=\nu u_{xx},\ in\ I\times (0,T],\\
(\rho,u)=(\rho_0,u_0), \ on \ I\times \{t=0\},\\
\rho=u=u_x=0,\ on \  \partial I\times (0,T],
\end{array}
\right.
\end{eqnarray}
where $\nu=2\mu+\lambda$.

The non-existence of Cauchy problem \eqref{1}-\eqref{3} in  $ C^1([0,T];H^m(\mathbb{R}))$, $m>2$ is  equivalent to the non-existence of  the  initial-boundary value problem \eqref{10} in  $C^{2,1}(\bar{I}\times [0,T])$,  which denotes the collection of functions that are $C^2$ in space and $C^1$ in time in $\bar{I}\times [0,T]$ here and in the following sections. Thus, in order to prove Theorem \ref{1001}, one needs only to show the following:
\begin{Theorem}\label{1000} The  initial-boundary value problem \eqref{10}  has no solution $(\rho,u)$ in  $C^{2,1}(\bar{I}\times [0,T])$ for any positive time  $T$, if the initial data $(\rho_0,u_0)$ satisfy the condition \eqref{11} or \eqref{12}.
\end{Theorem}

Let $\eta(x,t)$ denote the position
of the gas particle starting from $x$ at time $t=0$ satisfying
\begin{equation}\label{16.5}
\left\{ \begin{array}{ll}
\eta_t(x,t)=u(\eta(x,t),t), \\
\eta(x,0)=x.
\end{array}
\right.
\end{equation}
$\varrho$ and $v$ are the Lagrangian density and  velocity given by
\begin{eqnarray*}
\left\{ \begin{array}{ll}
\varrho(x,t)=\rho(\eta(x,t),t),\\
v(x,t)=u(\eta(x,t),t).
\end{array}
\right.
\end{eqnarray*}
Then  the system \eqref{10} can be rewritten in the Lagrangian coordinates as
\begin{eqnarray}\label{17}
\left\{ \begin{array}{ll}
\varrho_t+{\f{\varrho v_x}{\eta_x}}=0, \ in \ I\times(0,T],\\
\eta_x\varrho v_t+(\varrho^\gamma)_x=\nu({\f{v_x}{\eta_x}})_x, \ in \  I\times(0,T],\\
\eta_t(x,t)=v(x,t), \\
(\varrho,v,\eta)=(\rho_0,u_0,x), \ on \  I\times \{t=0\},\\
\varrho=v=v_x=0,\ on \  \partial I\times (0,T].
\end{array}
\right.
\end{eqnarray}
The first equation of  \eqref{17} implies that
\begin{eqnarray*}
\varrho(x,t)={\f{\rho_0(x)}{\eta_x(x,t)}}.
\end{eqnarray*}
Regarding $\rho_0$  as a parameter,  then one can reduce the system $\eqref{17}$  further  to
\begin{equation}\label{18}
\left\{ \begin{aligned}
&\rho_0v_t+(\f{\rho_0^\gamma}{\eta_x^\gamma})_x=\nu(\f{v_x}{\eta_x})_x, \ in \ I\times(0,T],\\
&\eta_t(x,t)=v(x,t), \\
& (v,\eta)=(u_0,x), \ on \  I\times \{t=0\},\\
&v=v_x=0,\ on \  \partial I\times (0,T].
\end{aligned}
\right.
\end{equation}

The condition \eqref{11} or \eqref{12} on the initial data $(\rho_0,u_0)$ takes the following form in the Lagrangian coordinates
\begin{eqnarray}\label{19}
\left\{ \begin{array}{ll}
\frac{(\rho_0)_x}{\rho_0}\geq \lambda_1,\ in \ (0,\lambda_3),\\
v_0(\lambda_3)<0,v_0\leq0,\ in \ (0,\lambda_3),
\end{array}
\right.
\end{eqnarray}
or
\begin{eqnarray}\label{20}
\left\{ \begin{array}{ll}
\frac{(\rho_0)_x}{\rho_0}\leq-\lambda_2,\ in \ (\lambda_4,1),\\
v_0(\lambda_4)>0,v_0\geq0,\ in \ (\lambda_4,1).
\end{array}
\right.
\end{eqnarray}

 The non-existence of the initial-boundary value problem \eqref{10}   is  equivalent to the non-existence of  the  initial-boundary value problem \eqref{18} in   $C^{2,1}(\bar{I}\times [0,T])$. Thus, Theorem \ref{1000} is a consequence of the following:

 \begin{Theorem}\label{1006} The problem \eqref{18}  has no solution $(v,\eta)$ in  $C^{2,1}(\bar{I}\times [0,T])$ for any positive time $T$, if the initial data $(\rho_0,u_0)$ satisfy the condition \eqref{19} or \eqref{20}.
\end{Theorem}

\subsection{Proof of Theorem \ref{1006}}
Given a sufficiently small positive time $T^*$, we let $(v,\eta)\in C^{2,1}(\bar{I}\times [0,T^*])$ be a solution of the system \eqref{18} with  \eqref{19} or \eqref{20}.
Define the  linear parabolic operator  $\rho_0\partial_t+L$ by
\begin{eqnarray}\label{21}
\rho_0\partial_t+L:=\rho_0\partial_t-\frac{\nu}{\eta_x}\partial_{xx}+\frac{\nu \eta_{xx}}{\eta_x^2}\partial_x,
\end{eqnarray}
where
\begin{eqnarray*}
\eta_x=1+\int_0^tv_xds \quad\ \ and\ \quad \eta_{xx}=\int_0^tv_{xx}d s.
\end{eqnarray*}
Then,  it  follows from the first equation of \eqref{18} that
\begin{eqnarray}\label{22}
\rho_0v_t+Lv=-(\frac{\rho_0^\gamma}{\eta_x^\gamma})_x.
\end{eqnarray}
Let $M$ be a positive constant such that
\begin{eqnarray}\label{23}
\rho_0+|v_0|+|(v_0)_x|+|(v_0)_{xx}|<M.
\end{eqnarray}
It follows from the continuity on time that for short time, it holds that
\begin{eqnarray*}
|v|+|v_x|+|v_{xx}|\leq M,\ in \ I\times (0,T^*].
\end{eqnarray*}
Taking a positive time $T<T^*$  sufficiently small such that $T\leq\frac{1}{2M}$, then one has
\begin{eqnarray}\label{24}
|\int_0^tv_xds|\leq MT\leq \frac{1}{2},\ in \ I\times (0,T].
\end{eqnarray}
This implies
\begin{eqnarray*}
\frac{1}{2}\leq\eta_x\leq\frac{3}{2},\ in \ I\times (0,T].
\end{eqnarray*}
Thus, the equation \eqref{22}  is a well-defined integro-differential equation  with degeneracy  for t-derivative due to that the initial density $\rho_0$ vanishes on the boundary $\partial I$.

Restrict $T$ further  such that $T\leq \frac{\lambda_1}{4M}$. Then,    \eqref{24} implies
\begin{eqnarray}\label{25}
-(\frac{\rho_0^\gamma}{\eta_x^\gamma})_x
=-\frac{\gamma\rho_0^\gamma}{\eta_x^\gamma}[\frac{(\rho_0)_x}{\rho_0}-\frac{\eta_{xx}}{\eta_x}]
\leq-\frac{\gamma\rho_0^\gamma}{\eta_x^\gamma}(\lambda_1-\frac{\lambda_1}{2})<0, \ in \ (0,\lambda_3)\times (0,T].
\end{eqnarray}
Thus, it follows from \eqref{22} and \eqref{25} that $v$ satisfies the following differential inequality
\begin{eqnarray}\label{26}
\rho_0v_t+Lv
\leq0, \ in \ (0,\lambda_3)\times (0,T].
\end{eqnarray}
Similarly, $v$ also satisfies
\begin{eqnarray}\label{27}
\rho_0v_t+Lv
\geq0, \ in \ (\lambda_4,1)\times (0,T].
\end{eqnarray}

In the rest of this section, our main task is to establish the Hopf's lemma and a strong maximum principle for the  differential inequality \eqref{26} and \eqref{27}.
First recall the definition of the parabolic boundary (see \cite{H}) of a bounded domain $D$ of $\mathbb{R}^d\times \mathbb{R}^+$. The parabolic boundary $\partial_p D$ of $D$ consists of  points $(x_0,t_0)\in \partial D$ such that $B_r(x_0)\times (t_0-r^2,t_0]$
contains points not in $D$, for any $r>0$.   In the following, suppose that $U$  is a bounded domain of $\mathbb{R}$, we use the notation $U_T:=U\times (0,T]$ to denote the cylinder in $(0,\lambda_3)\times (0,T]$.
Let $Q_T$ be any domain contained in $(0,\lambda_3)\times (0,T]$. We then derive a weak maximum principle for the  differential inequality \eqref{26} in $Q_T$.
\begin{Lemma}\label{1007}
Suppose that $w\in C^{2,1}(Q_T)\cap C(\bar{Q}_T)$ satisfies
\begin{eqnarray}\label{28}
\rho_0w_t+Lw
\leq0, \  in \ Q_T.
\end{eqnarray}
Then $w$ attains its  maximum on the parabolic boundary of $Q_T$.
\end{Lemma}
\textbf{Proof.} We first prove the statement under a stronger hypothesis instead of \eqref{28} that
\begin{eqnarray}\label{29}
\rho_0w_t+Lw<0, \  in \ Q_T.
\end{eqnarray}
Assume $w$ attains its  maximum at an interior point $(x_0,t_0)$ of the domain $Q_T$. Therefore
\begin{eqnarray*}
w_t(x_0,t_0)\geq0,  w_x(x_0,t_0)=0, w_{xx}(x_0,t_0)\leq0,
\end{eqnarray*}
which implies $\rho_0w_t+Lw\geq0$, this contradicts  \eqref{29}. Next, define the auxiliary function
\begin{eqnarray*}
\varphi^\varepsilon=w-\varepsilon t,
\end{eqnarray*}
for a  positive number $\varepsilon$. Then
\begin{eqnarray*}
\rho_0\varphi^\varepsilon_t+L\varphi^\varepsilon&=&\rho_0w_t+Lw-\varepsilon \rho_0<0, \  in \ Q_T.
\end{eqnarray*}
Thus $\varphi^\varepsilon$ attains its  maximum on the parabolic boundary of $Q_T$, which proves the assertion of Lemma \ref{1007} by letting $\varepsilon$ go to zero.\qed

The result in Lemma \ref{1007} can be extended to a general   domain $D\subset (0,\lambda_3)\times (0,T]$ (see \cite{Friedman}).
\begin{Lemma}\label{1008}
Suppose that $w\in C^{2,1}(D)\cap C(\bar{D})$ satisfies
\begin{eqnarray}\label{30}
\rho_0w_t+Lw
\leq0, \  in \ D.
\end{eqnarray}
Then $w$ attains its  maximum on the parabolic boundary  of $D$.
\end{Lemma}

 Next, we prove the  Hopf's  lemma for the  differential inequality \eqref{26}, which is critical for proving Theorem \ref{1006}.
\begin{Proposition}\label{1009}
Suppose that $w\in C^{2,1}((0,\lambda_3)\times (0,T])\cap C([0,\lambda_3]\times [0,T])$ satisfies \eqref{26} and there exits a point $(0,t_0)\in\{0\}\times (0,T]$ such that $w(x,t)<w(0,t_0)$ for any point $(x,t)$ in a
neighborhood $D$ of the point $(0,t_0)$, where
\begin{eqnarray*}
D=:\{(x,t): (x-r)^2+(t_0-t)< r^2,0< x<\frac{r}{2},0<t\leq t_0\}, 0<r<\lambda_3.
\end{eqnarray*}
 Then it holds that
\begin{eqnarray*}
\frac{\partial w(0,t_0)}{\partial \vec{n}}>0,
\end{eqnarray*}
where $\vec{n}$ is the outer unit normal vector at the point $(0,t_0)$.
\end{Proposition}
\textbf{Proof.}
For positive constants $\alpha$ and $\varepsilon$ to be determined,  set
\begin{eqnarray*}
q(\alpha,x,t)=e^{-\alpha[(x-r)^2+(t_0-t)]}-e^{-\alpha r^2}
\end{eqnarray*}
and
\begin{eqnarray*}
\varphi(\varepsilon,\alpha,x,t)=w(x,t)-w(0,t_0)+\varepsilon q(\alpha,x,t).
\end{eqnarray*}

First, we determine  $\varepsilon$.
The parabolic boundary $\partial_p D$ consists of two parts $\Sigma_1$ and $\Sigma_2$ given by
\begin{gather*}
\Sigma_1=\{(x,t):(x-r)^2+(t_0-t)<r^2, x=\frac{r}{2}, 0<t\leq t_0\}
\end{gather*}
and
\begin{gather*}
\Sigma_2=\{(x,t):(x-r)^2+(t_0-t)=r^2, 0\leq x\leq\frac{r}{2}, 0<t\leq t_0\}.
\end{gather*}
On  $\bar{\Sigma}_1$, $w(x,t)-w(0,t_0)<0$, and hence $w(x,t)-w(0,t_0)<-\varepsilon_0$ for some $\varepsilon_0>0$.
Note that $q\leq1$ on $\Sigma_1$. Then for such an $\varepsilon_0$,  $\varphi(\varepsilon_0,\alpha,x,t)<0$ on $\Sigma_1$. For $(x,t)\in \Sigma_2$, $q=0$ and $w(x,t)\leq w(0,t_0)$. Thus, $\varphi(\varepsilon_0,\alpha,x,t)\leq0$
for any $(x,t)\in \Sigma_2$ and $\varphi(\varepsilon_0,\alpha,0,t_0)=0$. One concludes that
\begin{eqnarray}\label{31}
\left\{ \begin{array}{ll}
\varphi(\varepsilon_0,\alpha,x,t)\leq0,\  on\  \partial_p D,\\
\varphi(\varepsilon_0,\alpha,0,t_0)=0.
\end{array}
\right.
\end{eqnarray}

Next, we choose  $\alpha$. It follows from \eqref{26} that
\begin{eqnarray}\label{32}
&&\rho_0\varphi_t(\varepsilon_0,\alpha,x,t)+L\varphi(\varepsilon_0,\alpha,x,t)\nonumber\\
&=&\rho_0w_t(x,t)+L w(x,t)
+\varepsilon_0[\rho_0q_t(\alpha,x,t)+L q(\alpha,x,t)]\nonumber\\
&\leq&\varepsilon_0[\rho_0q_t(\alpha,x,t)+L q(\alpha,x,t)].
\end{eqnarray}
A direct calculation yields
\begin{eqnarray}\label{33}
&&e^{\alpha[(x-r)^2+(t_0-t)]}[\rho_0q_t(\alpha,x,t)+L q(\alpha,x,t)]\nonumber\\
&=&-\frac{4\nu(x-r)^2}{\eta_x}\alpha^2
+[\rho_0+\frac{2\nu}{\eta_x}+\frac{2\nu \eta_{xx} (r-x)}{\eta_x^2}]\alpha\nonumber\\
&\leq& -\frac{2\nu r^2}{3}\alpha^2
+(M+4\nu+8\nu M r)\alpha.
\end{eqnarray}
Therefore, there exists a positive number $\alpha_0=\alpha_0(\nu,r,M)$ such that
\begin{eqnarray}\label{34}
\rho_0q_t(\alpha_0,x,t)+L q(\alpha_0,x,t)\leq0,\ in \ D,
\end{eqnarray}
Thus, it follows  from  \eqref{32} and \eqref{34} that
\begin{eqnarray}\label{35}
\rho_0\varphi_t(\varepsilon_0,\alpha_0,x,t)+L\varphi(\varepsilon_0,\alpha_0,x,t)\leq 0,\ in \ D.
\end{eqnarray}
 In conclusion, in view of  \eqref{31} and \eqref{35}, one has
\begin{eqnarray*}
\left\{ \begin{array}{ll}
\rho_0\varphi_t(\varepsilon_0,\alpha_0,x,t)+L\varphi(\varepsilon_0,\alpha_0,x,t)\leq 0,\ in \ D,\\
\varphi(\varepsilon_0,\alpha_0,x,t)\leq0, \ on \ \partial_pD,\\
\varphi(\varepsilon_0,\alpha_0,0,t_0)=0.
\end{array}
\right.
\end{eqnarray*}
This, together with Lemma \ref{1008} yields
\begin{eqnarray*}
\varphi(\varepsilon_0,\alpha_0,x,t)\leq0, \ in\ D.
\end{eqnarray*}
Therefore,  $\varphi(\varepsilon_0,\alpha_0,\cdot,\cdot)$ attains its maximum at the point $(0,t_0)$ in $D$. In particular, it holds that
\begin{eqnarray*}
\varphi(\varepsilon_0,\alpha_0,x,t_0)\leq \varphi(\varepsilon_0,\alpha_0,0,t_0)\quad \ for\ all\  x\in (0,\frac{r}{2}).
\end{eqnarray*}
This  implies
\begin{eqnarray*}
\frac{\partial \varphi(\varepsilon_0,\alpha_0,0,t_0)}{\partial \vec{n}}\geq0.
\end{eqnarray*}
Finally,  we get
\begin{eqnarray*}
\frac{\partial w(0,t_0)}{\partial \vec{n}}\geq-\varepsilon_0 \frac{\partial q(\alpha_0,0,t_0)}{\partial \vec{n}}
=2\varepsilon_0\alpha_0 r e^{-\alpha_0 r^2}>0.
\end{eqnarray*}
\qed

In order to establish a strong maximum principle for the  differential inequality \eqref{26}, we need to study the t-derivative of
interior maximum point. The main ideas in the following lemmas come from \cite{Friedman}.
\begin{Lemma}\label{1010}
Let  $w\in C^{2,1}((0,\lambda_3)\times (0,T])\cap C([0,\lambda_3]\times [0,T])$ satisfy \eqref{26} and have a maximum $M_0$ in the domain $(0,\lambda_3)\times (0,T]$. Suppose that $(0,\lambda_3)\times (0,T]$ contains a closed solid ellipsoid
\begin{eqnarray*}
\Omega^\sigma:=\{(x,t): (x-x_*)^2+\sigma(t-t_*)^2\leq r^2\},\sigma>0
\end{eqnarray*}
and  $w(x,t)<M_0$ for any interior point $(x,t)$ of $\Omega^\sigma$ and  $w(\bar{x},\bar{t})=M_0$ at some point $(\bar{x},\bar{t})$ on the  boundary   of $\Omega^\sigma$. Then $\bar{x}=x_*$.
\end{Lemma}
\textbf{Proof.}  Without loss of generality, one can assume that $(\bar{x},\bar{t})$ is the only point on $\partial \Omega^\sigma$ such that $w=M_0$ in $\Omega^\sigma$.  Otherwise,  one can limit it to a smaller closed ellipsoid lying in $\Omega^\sigma$ and having $(\bar{x},\bar{t})$ as the only common point with $\partial \Omega^\sigma$. We prove the desired result by contradiction.
Suppose that $\bar{x}\neq x_*$. Applying Lemma \ref{1008} on $\Omega^\sigma$ shows $\bar{t}<T$.
 Choose a  closed ball $D$ with center $(\bar{x},\bar{t})$ and radius $\tilde{r}<|\bar{x}-x_*|$ contained in $(0,\lambda_3)\times (0,T]$. Then $|x-x_*|\geq|\bar{x}-x_*|-\tilde{r}=:\hat{r}$ for any point $(x,t)\in D$.
The parabolic boundary  of $D$ is composed of a part $\Sigma_1$ lying in $\Omega^\sigma$ and a part $\Sigma_2$ lying outside $\Omega^\sigma$.

For positive constants $\alpha$ and $\varepsilon$ to be determined, set
\begin{eqnarray*}
q(\alpha,x,t)=e^{-\alpha[(x-x_*)^2+\sigma(t-t_*)^2]}-e^{-\alpha r^2}
\end{eqnarray*}
and
\begin{eqnarray*}
\varphi(\varepsilon,\alpha,x,t)=w(x,t)-M_0+\varepsilon q(\alpha,x,t).
\end{eqnarray*}

Note that $q(\alpha,x,t)>0$ in the interior of $\Omega^\sigma$, $q(\alpha,x,t)=0$ on $\partial \Omega^\sigma$ and $q(\alpha,x,t)<0$ outside $\Omega^\sigma$. So, it holds that $\varphi(\varepsilon,\alpha,\bar{x},\bar{t})=0$.
On  $\Sigma_1$, $w(x,t)-M_0<0$, and hence $w(x,t)-M_0<-\varepsilon_0$ for some $\varepsilon_0>0$.
Note that $q(\alpha,x,t)\leq1$ on $\Sigma_1$. Then for such an $\varepsilon_0$,  $\varphi(\varepsilon_0,\alpha,x,t)<0$ on $\Sigma_1$. For $(x,t)\in \Sigma_2$,  $q(\alpha,x,t)<0$ and $w(x,t)-M_0\leq0$. Thus, $\varphi(\varepsilon_0,\alpha,x,t)<0$
for any $(x,t)\in \Sigma_2$. One concludes that
\begin{eqnarray}\label{36}
\left\{ \begin{array}{ll}
\varphi(\varepsilon_0,\alpha,x,t)<0,\  on\  \partial_pD,\\
\varphi(\varepsilon_0,\alpha,\bar{x},\bar{t})=0.
\end{array}
\right.
\end{eqnarray}

Next, we estimate $\rho_0q_t(\alpha,x,t)+L q(\alpha,x,t)$.
One calculates that for $(x,t)\in D$,
\begin{eqnarray*}
&&e^{\alpha[(x-x_*)^2+\sigma(t-t_*)^2]}[\rho_0q_t(\alpha,x,t)+L q(\alpha,x,t)]\nonumber\\
&=&-\frac{4\nu(x-x_*)^2}{\eta_x}\alpha^2
+[2\sigma\rho_0(t_*-t)+\frac{2\nu}{\eta_x}+\frac{2\nu \eta_{xx} (x_*-x)}{\eta_x^2}]\alpha\nonumber\\
&\leq& -\frac{8\nu \hat{r}^2}{3}\alpha^2
+(2\sigma M+4\nu+8\nu M r)\alpha.
\end{eqnarray*}
Therefore, there exists a positive number $\alpha_0=\alpha_0(\nu,r,\hat{r},\sigma,M)$ such that
\begin{eqnarray}\label{37}
\rho_0q_t(\alpha_0,x,t)+L q(\alpha_0,x,t)\leq0,\ in \ D.
\end{eqnarray}
Thus, it follows  from  \eqref{26}, \eqref{32} and \eqref{37} that
\begin{eqnarray}\label{38}
\rho_0\varphi_t(\varepsilon_0,\alpha_0,x,t)+L\varphi(\varepsilon_0,\alpha_0,x,t)\leq 0,\ in \ D.
\end{eqnarray}
 In conclusion, it follows  from \eqref{36} and \eqref{38} that
\begin{eqnarray*}
\left\{ \begin{array}{ll}
\rho_0\varphi_t(\varepsilon_0,\alpha_0,x,t)+L\varphi(\varepsilon_0,\alpha_0,x,t)\leq 0,\ in \ D,\\
\varphi(\varepsilon_0,\alpha_0,x,t)<0, \ on \ \partial_pD,\\
\varphi(\varepsilon_0,\alpha_0,\bar{x},\bar{t})=0.
\end{array}
\right.
\end{eqnarray*}
However, Lemma \ref{1008} implies that
\begin{eqnarray*}
\varphi(\varepsilon_0,\alpha_0,x,t)<0, \ in\ D,
\end{eqnarray*}
which contradicts to $\varphi(\varepsilon_0,\alpha_0,\bar{x},\bar{t})=0$ due to $(\bar{x},\bar{t})\in D$ .\qed

Based on Lemma \ref{1010}, it is standard to prove the following lemma.  For details,  please refer to Lemma 3 of Chapter 2 in \cite{Friedman}.
\begin{Lemma}\label{1011}
Suppose that $w\in C^{2,1}((0,\lambda_3)\times (0,T])\cap C([0,\lambda_3]\times [0,T])$ satisfies \eqref{26}.
If $w$ has a maximum in an interior point $P_0=(x_0,t_0)$ of $(0,\lambda_3)\times (0,T]$, then
$w(P)=w(P_0)$ for any point $P(x,t_0)$ of $(0,\lambda_3)\times (0,T]$.
\end{Lemma}

We first prove a localized version strong maximum principle in a rectangle $\mathcal{R}$ of  the domain $(0,\lambda_3)\times (0,T]$.
\begin{Lemma}\label{1012}
Suppose that $w\in C^{2,1}((0,\lambda_3)\times (0,T])\cap C([0,\lambda_3]\times [0,T])$ satisfies \eqref{26}.
If $v$ has a maximum in the interior point $P_0=(x_0,t_0)$ of $(0,\lambda_3)\times (0,T]$, then there exists a rectangle
\begin{eqnarray*}
\mathcal{R}(P_0):=\{(x,t):x_0-a_1\leq x\leq x_0+a_1, t_0-a_0\leq t\leq t_0\}
\end{eqnarray*}
 in $(0,\lambda_3)\times (0,T]$ such that
$w(P)=w(P_0)$ for any point $P$ of $\mathcal{R}(P_0)$.
\end{Lemma}
\textrm{\textbf{Proof}} We prove the desired result by contradiction. Suppose that there exists an interior point $P_1=(x_1,t_1)$ of $(0,\lambda_3)\times (0,T]$ with $t_1< t_0$ such that $w(P_1)<w(P_0)$.
 Connect  $P_1$ to $P_0$ by a simple smooth curve $\gamma$. Then there exists a point $P_*=(x_*,t_*)$ on $\gamma$ such that $w(P_*)=w(P_0)$ and $w(\bar{P})<w(P_*)$ for all any point $\bar{P}$ of $\gamma$ between $P_1$ and $P_*$. We may assume that $P_*=P_0$ and $P_1$ is very near to $P_0$. There exist a rectangle $\mathcal{R}(P_0)$ in $(0,\lambda_3)\times (0,T]$ with small positive numbers $a_0$ and $a_1$ (will be determined) such that
 $P_1$ lies on $t=t_0-a_0$.
  Since $\mathcal{R}(P_0)\setminus \{t=t_0\}\cap\{t=\bar{t}\}$ contains some point $\bar{P}=(\bar{x},\bar{t})$ of $\gamma$ and $w(\bar{P})<w(P_0)$,  we deduce $w(P)<w(P_0)$ for each point $P$ in $\mathcal{R}(P_0)\setminus \{t=t_0\}\cap\{t=\bar{t}\}$ due to Lemma  \ref{1011}.  Therefore, $w(P)<w(P_0)$ for each point $P$ in $\mathcal{R}(P_0)\setminus \{t=t_0\}$.

For positive constants $\alpha$ and $\varepsilon$ to be determined, set
\begin{eqnarray*}
q(\alpha,x,t)=t_0-t-\alpha(x-x_0)^2
\end{eqnarray*}
and
\begin{eqnarray*}
\varphi(\varepsilon,\alpha,x,t)=w(x,t)-w(P_0)+\varepsilon q(\alpha,x,t).
\end{eqnarray*}
Assume further that $P=(x_0-a_1,t_0-a_0)$ is on the parabola $q(\alpha,x,t)=0$. Then
\begin{equation}\label{38.2}
  \alpha=\frac{a_0}{a_1^2}.
\end{equation}

To choose  $\alpha$, one calculates
\begin{eqnarray}\label{38.4}
&&\rho_0q_t(\alpha,x,t)+L q(\alpha,x,t)\nonumber\\
&=&
-\rho_0+[\frac{2\nu}{\eta_x}-\frac{2\nu \eta_{xx} (x-x_0)}{\eta_x^2}]\alpha\nonumber\\
&\leq& -\rho_0
+(4\nu+8\nu M a_1)\alpha.
\end{eqnarray}
since $\rho_0$ has a positive lower bound depending on $x_0-a_1$ in $\mathcal{R}(P_0)$, one can choose $\alpha_0$  such that
\begin{equation}\label{39}
\alpha_0<\frac{\rho_0}{4\nu+8\nu M a_1}.
\end{equation}
This and \eqref{38.4} imply that
\begin{eqnarray}\label{40}
\rho_0\varphi_t(\alpha_0,x,t)+L \varphi(\alpha_0,x,t)
\leq0, \ in\ \mathcal{R}(P_0).
\end{eqnarray}

One can now fix $a_1$ such that
\begin{equation*}
  a_1<\min\{x_0,\lambda_3-x_0\}
\end{equation*}
and it then follows from \eqref{38.2} and \eqref{38.4} that one can choose $a_0$ such that
\begin{equation*}
  a_0<\min\{t_0,\frac{a_1^2\rho_0}{2(4\nu+8\nu M a_1)}\}.
\end{equation*}

Denote $\mathcal{S}=\{(x,t)\in \mathcal{R}(P_0), q(\alpha_0,x,t)\geq0\}$.  The parabolic boundary $\partial_p\mathcal{S}$ of $\mathcal{S}$ is composed of a part $\Sigma_1$ lying in $\mathcal{R}(P_0)$ and a part $\Sigma_2$ lying on  $\mathcal{R}(P_0)\cap \{t=t_0-a_0\}$.

We now determine $\varepsilon$. Note that on  $\Sigma_2$, $w(x,t)-M_0<0$, and $q(\alpha_0,x,t)$ is bounded, one can choose sufficiently small number $\varepsilon_0$ such that $\varphi(\varepsilon_0,\alpha_0,x,t)<0$ on $\Sigma_2$.
{On $\Sigma_1\setminus \{P_0\}$, $q(\alpha_0,x,t)=0$ and $w(x,t)-M_0<0$}. Thus, $\varphi(\varepsilon_0,\alpha_0,x,t)<0$ on $\Sigma_1\setminus \{P_0\}$ and
 $\varphi(\varepsilon_0,\alpha_0,x_0,t_0)=0$. One concludes that
\begin{eqnarray}\label{40.5}
\left\{ \begin{array}{ll}
\varphi(\varepsilon_0,\alpha_0,x,t)<0,\  on\  \partial_p\mathcal{S}\setminus \{P_0\},\\
 \varphi(\varepsilon_0,\alpha_0,x_0,t_0)=0.
\end{array}
\right.
\end{eqnarray}

 In conclusion,  it follows from \eqref{40} and \eqref{40.5} that there exist $\varepsilon_0$, $a_0$ and $a_1$ such that
\begin{eqnarray}\label{41}
\left\{ \begin{array}{ll}
\rho_0\varphi_t(\varepsilon_0,\alpha_0,x,t)+L\varphi(\varepsilon_0,\alpha_0,x,t)\leq 0,\ in\ \mathcal{S},\\
\varphi(\varepsilon_0,\alpha_0,x,t)<0,\  on\  \partial_p\mathcal{S}\setminus \{P_0\},\\
 \varphi(\varepsilon_0,\alpha_0,x_0,t_0)=0.
\end{array}
\right.
\end{eqnarray}

In view of Lemma \ref{1008} and \eqref{41},  $\varphi(\varepsilon_0,\alpha_0,\cdot,\cdot)$ only attains its maximum at $P_0$ in $\mathcal{\bar{S}}$, thus
\begin{eqnarray*}
\frac{\partial \varphi(\varepsilon_0,\alpha_0,x_0,t_0)}{\partial t}\geq0.
\end{eqnarray*}
Note that $q$ satisfies at $P_0$
\begin{eqnarray*}
\frac{\partial q(\alpha_0,x_0,t_0)}{\partial t}=-1.
\end{eqnarray*}
 Therefore
\begin{eqnarray}\label{42}
\frac{\partial w(x_0,t_0)}{\partial t}\geq \varepsilon_0.
\end{eqnarray}
But, by the assumption, $w$  attains its   maximum   at $P_0$,  it follows that
\begin{eqnarray*}
\rho_0\frac{\partial w(x_0,t_0)}{\partial t}\leq -L w(x_0,t_0)\leq0,
\end{eqnarray*}
which contradicts \eqref{42}.\qed

 Now we  can prove the following strong maximum principle.
\begin{Proposition}\label{1013}
Suppose that $w\in C^{2,1}((0,\lambda_3)\times (0,T])\cap C([0,\lambda_3]\times [0,T])$ satisfies \eqref{26}.
If $w$ attains its maximum at some interior point $P_0=(x_0,t_0)$ of $ (0,\lambda_3)\times (0,T]$, then
$w(P)=w(P_0)$ for any point  $P\in(0,\lambda_3)\times (0,t_0]$.
\end{Proposition}
\textrm{\textbf{Proof}} We prove the desired result by contradiction.  Suppose that $w\not\equiv w(P_0)$. Then there exists a point $P_1=(x_1,t_1)$ of $ (0,\lambda_3)\times (0,t_0]$ such that $w(P_1)<w(P_0)$. By Lemma \ref{1011}, there must be $t_1<t_0$.

Connect $P_1$ to $P_0$ by a straight line $\gamma$.  There exists a point $P_*$ on $\gamma$ such that $w(P_*)=w(P_0)$ and $w(\bar{P})<w(P_*)$ for any point $\bar{P}$ on $\gamma$ lying between $P_*$ and $P_1$. Denote by $\gamma_0$ the closed sub straight line of $\gamma$ lying $P_*$ and $P_1$.
Construct a series of rectangles $\mathcal{R}_n, n=1,2,\cdots,N$ with  small  $a_n$ and $b_n$ such that $\gamma_0\subset\cup_{n=1}^N\mathcal{R}_n$, $P_*\in \mathcal{R}_1$ and $P_1\in \mathcal{R}_N$. Applying Lemma \ref{1012} on  $\mathcal{R}_1, \mathcal{R}_2, \cdots, \mathcal{R}_N$ step by step it follows that $w=w(P_1)$ in  $\cup_{n=1}^N\mathcal{R}_n$. Hence, one deduces $w(P_*)\equiv w(P_1)$ due to $P_*$ lying on  $\gamma_0$, which is a contradiction.\qed

Let $D$ be a bounded domain contained in the domain $ (\lambda_4,1)\times (0,T]$. Similar to Lemma \ref{1008}, Proposition \ref{1009} and Proposition \ref{1013},
we have  corresponding  weak maximum principle,  Hopf's lemma and strong minimum principle for the  differential inequality \eqref{27}.
\begin{Lemma}\label{1014}
Suppose that $w\in C^{2,1}(D)\cap C(\bar{D})$ satisfies
\begin{eqnarray*}
\rho_0w_t+Lw
\geq0, \  in \ D.
\end{eqnarray*}
Then $w$ attains its  minimum on the parabolic boundary  of $D$.
\end{Lemma}

\begin{Proposition}\label{1015}
Suppose that $w\in C^{2,1}((\lambda_4,1)\times (0,T])\cap C([\lambda_4,1]\times [0,T])$ satisfies \eqref{27} and there exits a point $(1,t_0)\in\{1\}\times (0,T]$ such that $w(x,t)>w(1,t_0)$ for any point $(x,t)$ in a
neighborhood $D$ of the point $(0,t_0)$, where
\begin{eqnarray*}
D=:\{(x,y):(x-(1-r))^2+(t_0-t)< r^2, 1-\frac{r}{2}<x<1 , 0<t\leq t_0\}, 1-r>\lambda_4.
\end{eqnarray*}
Then it holds that
\begin{eqnarray*}
\frac{\partial w(1,t_0)}{\partial \vec{n}}<0,
\end{eqnarray*}
where $\vec{n}$ is the outer unit normal vector at the point $(1,t_0)$.
\end{Proposition}

\begin{Proposition}\label{1016}
Suppose that $w\in C^{2,1}((\lambda_4,1)\times (0,T])\cap C([\lambda_4,1]\times [0,T])$ satisfies \eqref{27}.
If $w$ attains its minimum at some interior point $P_0=(x_0,t_0)$ of $ (\lambda_4,1)\times (0,T]$, then
$w(P)=w(P_0)$ for any point  $P$ of $(\lambda_4,1)\times (0,t_0]$.
\end{Proposition}

We are now ready to prove Theorem \ref{1006}.\\
\textrm{\textbf{Proof of Theorem \ref{1006}.}} We first consider the case of the domain $(0,\lambda_3)\times (0,T]$.
 Recall $v$ satisfies \eqref{26}, so the weak maximum principle, Hopf lemma  and strong maximum principle for $w$ holds also for $v$.
Since  $v_0(\lambda_3)<0$, by continuity of $v$ on time, then there exists a time $t_0>0$ such that $v(\lambda_3,\cdot)<0$ in $ (0,t_0)$. By Lemma \ref{1007},  $v$  attains its  maximum on the parabolic boundary $\{x=0\}\times (0,t_0]\cup \{x=\lambda_3\}\times (0,t_0]\cup (0,\lambda_3)\times\{t=0\}$.  Since $v=0$ on the parabolic boundary $\{x=0\}\times (0,t_0]$ and $v_0\leq0$ in $ (0,\lambda_3)$, by Proposition \ref{1013},  $v$ only attains its  maximum on the set $\{x=0\}\times (0,t_0]\cup (0,\lambda_3)\times\{t=0\}$.  Thus, $v(x,t)<v(0,t_0)(=0)$ for any point $(x,t)\in (0,\lambda_3)\times (0,t_0]$.  Applying Proposition \ref{1009} shows that $\frac{\partial v(0,t_0)}{\partial \vec{n}}>0$, which contradicts to $v_x(x,t)=0$ on $\partial I\times (0,T]$ of the system \eqref{18}. The other case is similar.\qed

\begin{center}
\section{Proof of Theorem \ref{1003}}
\end{center}

\subsection{Reformulation of Theorem \ref{1003}}

Suppose that $\kappa=0$ and $n=1$.
Let $(\rho,u,e)\in C^1([0,T];H^m(\mathbb{R})), m>2$  be a solution  to the system \eqref{5}-\eqref{7} with  the initial density satisfying $\eqref{8}$.
Let $a(t)$ and $b(t)$ be the particle paths stating from $0$ and  $1$, respectively.
Similar to \eqref{9}, one can show that
\begin{eqnarray*}
\left\{ \begin{array}{ll}
[a(t), b(t)]=[0,1],\\
u(x,t)=u_x(x,t)=0.
\end{array}
\right.
\end{eqnarray*}
where $t\in(0,T^*)$ and  $x\in [a(t), b(t)]^c$.

Therefore, to study the ill-posedness of the system \eqref{5}-\eqref{7}  with  the initial density satisfying $\eqref{8}$ is  equivalent to
study that of the following initial-boundary value problem
\begin{eqnarray}\label{13}
\left\{ \begin{array}{ll}
\rho_t+(\rho
u)_x=0,\ in\ I\times (0,T],\\
(\rho u)_t+(\rho u^2+p)_x=\mu u_{xx},\ in\ I\times (0,T],\\
(\rho e)_t+(\rho eu)_x+ pu_x=\mu u_x^2,\ in\ I\times (0,T],\\
(\rho,u,e)=(\rho_0,u_0,e_0), \ on \ I\times \{t=0\},\\
\rho=u=u_x=0,\ on \  \partial I\times (0,T].
\end{array}
\right.
\end{eqnarray}

The non-existence of Cauchy problem \eqref{5}-\eqref{7} in   $ C^1([0,T];H^m(\mathbb{R}))$, $m>2$ is  equivalent to the non-existence of  the  initial-boundary value problem \eqref{13} in  $C^{2,1}(\bar{I}\times [0,T])$. Thus, in order to prove Theorem \ref{1003}, we need only to show the following:
\begin{Theorem}\label{1002} The initial-boundary value problem \eqref{13} has no solution $(\rho,u,e)$ in   $C^{2,1}(\bar{I}\times [0,T])$ for any positive time  $T$, if the initial data $(\rho_0,u_0,e_0)$ satisfy the condition \eqref{14} or \eqref{15}.
\end{Theorem}

Let $\eta(x,t)$ be the position
of the gas particle starting from $x$ at time $t=0$ defined by \eqref{16.5}.
Let $\varrho$, $v$ and $\mathfrak{e}$ be the Lagrangian density, velocity and internal energy, which are defined by
\begin{eqnarray}\label{42.5}
\left\{ \begin{aligned}
\varrho(x,t)=\rho(\eta(x,t),t),\\
v(x,t)=u(\eta(x,t),t),\\
\mathfrak{e}(x,t)=e(\eta(x,t),t).
\end{aligned}
\right.
\end{eqnarray}
Then  the system \eqref{13} can be rewritten in the Lagrangian coordinates as
\begin{equation}\label{43}
\left\{ \begin{aligned}
&\rho_0v_t+(\f{\rho_0\mathfrak{e}}{\eta_x})_x=\mu(\f{v_x}{\eta_x})_x, \ in \ I\times(0,T],\\
&\rho_0\mathfrak{e}_t+(\gamma-1)\frac{\rho_0\mathfrak{e}v_x}{\eta_x}=\mu\f{v_x^2}{\eta_x},\ in \ I \times (0,T],\\
&\eta_t(x,t)=v(x,t),\\
& (v,\mathfrak{e},\eta)=(u_0,e_0,x), \ on \  I\times \{t=0\},\\
&v=v_x=0,\ on \  \partial I\times (0,T].
\end{aligned}
\right.
\end{equation}

In the Lagrangian coordinates, the condition \eqref{14} or \eqref{15} on the initial data $(\rho_0,u_0,\mathfrak{e}_0)$ becomes
\begin{eqnarray}\label{44}
\left\{ \begin{array}{ll}
\frac{(\rho_0)_x}{\rho_0}+\frac{(\mathfrak{e}_0)_x}{\rho_0}\geq \lambda_5,\ in \ (0,\lambda_7),\\
v_0(\lambda_7)<0,v_0\leq0,\ in \ (0,\lambda_7),
\end{array}
\right.
\end{eqnarray}
or
\begin{eqnarray}\label{45}
\left\{ \begin{array}{ll}
\frac{(\rho_0)_x}{\rho_0}+\frac{(\mathfrak{e}_0)_x}{\rho_0}\leq-\lambda_6,\ in \ (\lambda_8,1),\\
v_0(\lambda_8)>0,v_0\geq0,\ in \ (\lambda_8,1),
\end{array}
\right.
\end{eqnarray}
respectively.

 The non-existence of the initial-boundary value  problem \eqref{43}   is  equivalent to the non-existence of  the  initial-boundary value problem \eqref{12} in  $C^{2,1}(\bar{I}\times [0,T])$. Thus, in order to prove Theorem \ref{1002}, we need only to show the following:
\begin{Theorem}\label{1017} The initial-boundary value problem \eqref{43} has no solution  $(v,\mathfrak{e},\eta)$ in  $C^{2,1}(\bar{I}\times [0,T])$ for any positive time $T$, if the initial data $(\rho_0,u_0)$ satisfy the condition \eqref{44} or \eqref{45}.
\end{Theorem}

\subsection{Proof of Theorem \ref{1017}}
 Given sufficiently small positive time $T^*$.  Let $(v,\mathfrak{e},\eta)\in C^{2,1}(\bar{I}\times [0,T^*])$ be a solution of the system \eqref{43} with  \eqref{44} or \eqref{45}.
Define the  linear parabolic operator  $\rho_0\partial_t+L$ similar to Subsection 3.1  by
\begin{eqnarray*}
\rho_0\partial_t+L:=\rho_0\partial_t-\frac{\mu}{\eta_x}\partial_{xx}+\frac{\mu \eta_{xx}}{\eta_x^2}\partial_x,
\end{eqnarray*}
Then,  it  follows from the first equation of \eqref{43} that
\begin{eqnarray}\label{47}
\rho_0v_t+Lv=-(\f{\rho_0\mathfrak{e}}{\eta_x})_x.
\end{eqnarray}
Let $M$ be a positive constant such that
\begin{eqnarray*}
\rho_0+|v_0|+|(v_0)_x|+|(v_0)_{xx}|+|\mathfrak{e}_0|+|(\mathfrak{e}_0)_x|<M.
\end{eqnarray*}
It follows from continuity on time that for suitably small $T^*$ that
\begin{eqnarray*}
|v|+|v_x|+|v_{xx}|++|\mathfrak{e}|+|\mathfrak{e}_x|\leq M,\ in \ I\times (0,T^*]
\end{eqnarray*}
and
\begin{eqnarray}\label{48}
\frac{(\rho_0)_x}{\rho_0}+\frac{\mathfrak{e}_x}{\rho_0}\geq \frac{\lambda_5}{2},\ in \ (0,\lambda_7)\times (0,T^*].
\end{eqnarray}
Taking a positive time $T<T^*$  sufficiently small such that $T\leq\frac{1}{2M}$, then one gets
\begin{eqnarray*}
|\int_0^tv_xds|\leq MT\leq \frac{1}{2},\ in \ I\times (0,T].
\end{eqnarray*}
This implies
\begin{eqnarray*}
\frac{1}{2}\leq\eta_x\leq\frac{3}{2},\ in \ I\times (0,T].
\end{eqnarray*}
Thus, \eqref{43}  are  well-defined integro-differential equations  with degeneracy  for t-derivative due to that the initial density $\rho_0$ vanishes on the boundary $\partial I$.

Take $T$  small further such that $T\leq \frac{\lambda_5}{8M}$. Therefore,  \eqref{48} implies
\begin{eqnarray}\label{49}
-(\f{\rho_0\mathfrak{e}}{\eta_x})_x
=-\frac{\rho_0\mathfrak{e}}{\eta_x}[\frac{(\rho_0)_x}{\rho_0}
+\frac{\mathfrak{e}_x}{\rho_0}-\frac{\eta_{xx}}{\eta_x}]
\leq-\frac{\rho_0\mathfrak{e}}{\eta_x}(\frac{\lambda_5}{2}-\frac{\lambda_5}{4})<0, \ in \ (0,\lambda_7)\times (0,T].
\end{eqnarray}
Thus, it follows from \eqref{47} and \eqref{49} that $v$ satisfies the following  differential inequality
\begin{eqnarray*}
\rho_0w_t+Lw
\leq0, \ in \ (0,\lambda_7)\times (0,T].
\end{eqnarray*}
Similarly, $v$ also satisfies
\begin{eqnarray*}
\rho_0w_t+Lw
\geq0, \ in \ (\lambda_8,1)\times (0,T].
\end{eqnarray*}

The rest is  the same as the proof of Theorem \ref{1006} in Subsection 2.2 and thus omitted.

\begin{center}
\section{Proof of Theorem \ref{1005}}
\end{center}

\subsection{Reformulation of Theorem \ref{1005}}

Suppose that $\kappa>0$.
Let $(\rho,u,e)\in C^1([0,T];H^m(\mathbb{R}^n)), m>[\frac{n}{2}]+2$  be a solution  to the system \eqref{5}-\eqref{7} with  the initial density satisfying \eqref{8}. Denote by $X(x_0,t)$ the particle
trajectory starting at $x_0$ when $t=0$, that is,
\begin{equation*}
\left\{ \begin{array}{ll}
\partial_tX(x_0,t)=u(X(x_0,t),t), \\
X(x_0,0)=x_0.
\end{array}
\right.
\end{equation*}
Set
\begin{equation*}
\Omega=\Omega(0)\quad\ and\ \quad \Omega(t)=\{x=X(x_0,t):x_0\in \Omega(0)\}.
\end{equation*}
It follows from the first equation of \eqref{5} that $\mbox{supp}_x\,\rho=\Omega(t)$. Under the assumption that the entropy $S(t,x)$ is finite in the vacuum domain $\Omega(t)^c$, then one deduces from the  equation of state \eqref{6} that
\begin{equation*}
e(x,t)=0\quad \ for \ x\in \Omega(t)^c.
\end{equation*}
Due to $e(\cdot,t)\in H^m(\mathbb{R}^n), m>[\frac{n}{2}]+2$,   one gets
\begin{eqnarray*}\label{49.2}
e_{x_i}(x,t)=e_{x_ix_j}(x,t)=0\quad  \ for \ x\in \Omega(t)^c, i,j=1,2,\cdots,n.
\end{eqnarray*}
It follows from the third equation of $\eqref{5}$ that
\begin{eqnarray}
\frac{\mu}{2}|\nabla u+\nabla u^T|^2+\lambda(\textrm{div}u)^2=0\quad  \ for \ x\in \Omega(t)^c.
\end{eqnarray}
Following the arguments in \cite{Xin}, one can calculate that
\begin{eqnarray}\label{49.4}
\frac{\mu}{2}|\nabla u+\nabla u^T|^2+\lambda(\textrm{div}u)^2
\geq
\left\{ \begin{array}{ll}
(2\mu+n\lambda)\sum_{i=1}^n(u_{x_i})^2+\mu\sum_{i> j}^n(u_{x_i}+u_{x_j})^2, \ if \ \lambda\leq0,\\
2\mu\sum_{i=1}^n(u_{x_i})^2+\mu\sum_{i> j}^n(u_{x_i}+u_{x_j})^2, \ if \ \lambda>0,
\end{array}
\right.
\end{eqnarray}
this, together with \eqref{49.2} implies
\begin{eqnarray*}
\partial_iu_j+\partial_ju_i=0\quad  \ for \ x\in \Omega(t)^c, i,j=1,2,\cdots,n.
\end{eqnarray*}
Because of $u(\cdot,t)\in H^m(\mathbb{R}^n), m>[\frac{n}{2}]+2$,  it holds that
\begin{eqnarray*}
u(x,t)=u_{x_i}(x,t)=u_{x_ix_j}(x,t)=0\quad  \ for \ x\in \Omega(t)^c, i,j=1,2,\cdots,n.
\end{eqnarray*}
Furthermore, one has $\Omega(t)=\Omega(0)$.

One concludes that
\begin{eqnarray*}
\left\{ \begin{array}{ll}
\Omega(t)=\Omega(0),\\
e(x,t)=e_{x_i}(x,t)=0,
\end{array}
\right.
\end{eqnarray*}
where $t\in(0,T^*)$ and  $x\in \Omega(t)^c, i=1,2,\cdots,n$.

Therefore,  to study the ill-posedness of the system \eqref{5}-\eqref{7}  with  the initial density satisfying \eqref{8}, one needs only to
study the ill-posedness of the following  initial-boundary value problem
\begin{eqnarray}\label{16}
\left\{ \begin{array}{ll}
\partial_t\rho+\textrm{div}(\rho
u)=0,\ in\ \Omega\times (0,T],\\
\partial_t(\rho u)+\textrm{div}(\rho u\otimes u)+\nabla p=\mu\Delta u+(\mu+\lambda)\nabla\textrm{div}u,\ in\ \Omega\times (0,T],\\
\partial_t(\rho e)+\textrm{div}(\rho eu)+ p\textrm{div}u=\frac{\mu}{2}|\nabla u+(\nabla u)^*|^2+\lambda(\textrm{div}u)^2\\
\quad\quad\quad\quad\quad\quad\quad\quad\quad\quad\quad\quad\quad\quad\quad\quad\quad\quad\quad+\frac{\kappa(\gamma-1) }{R}\Delta e, \ in\ \Omega\times (0,T],\\
(\rho,u,e)=(\rho_0,u_0,e_0), \ on \ \Omega\times \{t=0\},\\
e(x,t)=e_{x_i}(x,t)=0, \ on \  \partial \Omega\times (0,T].
\end{array}
\right.
\end{eqnarray}

The non-existence of Cauchy problem \eqref{5}-\eqref{7} in  $ C^1([0,T];H^m(\mathbb{R}^n))$, $m>[\frac{n}{2}]+2$ will follow from the non-existence of  the  initial-boundary value problem \eqref{16} in  $C^{2,1}(\bar{I}\times [0,T])$. Thus, in order to prove Theorem \ref{1005}, we need only to show the following theorem:
\begin{Theorem}\label{1004} The  initial-boundary value problem \eqref{16} in the case of $\kappa>0$ has no solution $(\rho,u,e)$ in  $C^{2,1}(\bar{\Omega}\times [0,T])$ for any  positive time $T$.
\end{Theorem}

 Let $\eta(x,t)$ denote the position
of the gas particle starting from $x$ at time $t=0$ defined by \eqref{16.5}.
Let $\varrho$, $v$ and $\mathfrak{e}$ be the Lagrangian density, velocity and internal energy, respectively, which are defined by \eqref{42.5}.
We will also use the following notations (see also \cite{CS,CLS,JM1,JM2})
\begin{eqnarray*}
\left\{ \begin{array}{ll}
J=\det D\eta  \quad (Jacobian\  determinant),\\
B=[D\eta]^{-1}  \quad (inverse\ of \ deformation\ tensor),\\
b=JB   \quad (transpose\ of \ cofactor\ matrix).
\end{array}
\right.
\end{eqnarray*}

We will always use the convention in this section that repeated Latin indices $i,j,k,$ etc., are summed from $1$ to $n$.
Then  the system \eqref{5} can be rewritten in the Lagrangian coordinates as
\begin{eqnarray}\label{51}
\begin{cases}
\partial_t\varrho+\varrho B_i^j\partial_jv^i=0, \ in \ \Omega \times (0,T],\\
\varrho \partial_tv^i+(\gamma-1)B_i^j\partial_j(\varrho \mathfrak{e})=\mu B_l^k\partial_k(B_l^j\partial_jv^i)+(\mu+\lambda)B_i^k\partial_k(B_l^j\partial_jv^l),\ in \ \Omega \times (0,T],\\
\varrho \partial_t\mathfrak{e}+(\gamma-1)\varrho \mathfrak{e} B_i^j\partial_jv^i
=\frac{\mu}{2}|B_l^j\partial_jv^i+(B_l^j\partial_jv^i)^*|^2+\lambda(B_i^j\partial_jv^i)^2\\
\quad\quad\quad\quad\quad\quad\quad\quad\quad\quad\quad\quad\quad\quad\quad\quad\quad\quad\quad+\frac{\kappa(\gamma-1) }{R}B_l^k\partial_k(B_l^j\partial_j\mathfrak{e}),\ in \ \Omega \times (0,T],\\
\eta_t(x,t)=v(x,t),\\
(\varrho,v,\mathfrak{e},\eta)=(\rho_0,u_0,e_0,x), \ on \  \Omega\times \{t=0\},\\
\mathfrak{e}(x,t)=\mathfrak{e}_{x_i}(x,t)=0, \ on \  \partial \Omega\times (0,T].
\end{cases}
\end{eqnarray}
 It follows from  \eqref{51} that
\begin{eqnarray*}
\varrho(x,t)={\f{\rho_0(x)}{J(x,t)}}.
\end{eqnarray*}
Regarding the initial density $\rho_0$  as a parameter, one can rewrite the system $\eqref{51}$ as
\begin{eqnarray}\label{52}
\left\{ \begin{array}{ll}
\rho_0\partial_tv^i+(\gamma-1)b_i^j\partial_j(J^{-1}\rho_0 \mathfrak{e} )=\mu b_l^k\partial_k(J^{-1}b_l^j\partial_jv^i)\\
\quad\quad\quad\quad\quad\quad\quad\quad\quad\quad\quad\quad\quad\quad\quad\quad+(\mu+\lambda)b_i^k\partial_k(J^{-1}b_l^j\partial_jv^l), \ in \ \Omega\times(0,T],\\
\rho_0\partial_t\mathfrak{e}+(\gamma-1)J^{-1}\rho_0 \mathfrak{e} b_i^j\partial_jv^i=\frac{\mu}{2}J^{-1}|b_l^j\partial_jv^i+(b_l^j\partial_jv^i)^*|^2+\lambda J^{-1}(b_i^j\partial_jv^i)^2\\
\quad\quad\quad\quad\quad\quad\quad\quad\quad\quad\quad\quad\quad\quad\quad\quad\quad+\frac{\kappa(\gamma-1) }{R}b_l^k\partial_k(J^{-1}b_l^j\partial_j\mathfrak{e}),\ in \ \Omega \times (0,T],\\
\eta_t(x,t)=v(x,t),\\
 (v,\mathfrak{e},\eta)=(u_0,e_0,x), \ on \  \Omega\times \{t=0\},\\
\mathfrak{e}(x,t)=\mathfrak{e}_{x_i}(x,t)=0, \ on \  \partial \Omega\times (0,T].
\end{array}
\right.
\end{eqnarray}

 The non-existence of the initial-boundary value problem \eqref{16} will be a consequence of the non-existence of  the   initial-boundary value problem \eqref{51} in  $C^{2,1}(\bar{\Omega}\times [0,T])$. Thus, in order to prove Theorem \ref{1004}, we need only to show the following:
\begin{Theorem}\label{1018} The  problem \eqref{52} in the case of $\kappa>0$ has no solution  $(v, \mathfrak{e},\eta)$ in $C^{2,1}(\bar{\Omega}\times [0,T])$ for any positive time $T$.
\end{Theorem}

\subsection{Proof of Theorem \ref{1018}}
Let $T^*$ be a given suitably small positive time. Let $(v,\mathfrak{e},\eta)\in C^{2,1}(\bar{\Omega}\times [0,T^*])$ be a solution of the system \eqref{52}.
Let $M$ be a positive constant such that
\begin{eqnarray*}
\rho_0+\sum_{|\alpha|\leq2}|D^\alpha v_0|+\sum_{|\alpha|\leq2}|D^\alpha \mathfrak{e}_0|< M.
\end{eqnarray*}
It follows from continuity on time that for short time $T^*$
\begin{eqnarray*}
\sum_{|\alpha|\leq2}|D^\alpha v|+\sum_{|\alpha|\leq2}|D^\alpha \mathfrak{e}|\leq M,\ in \ \Omega\times (0,T^*].
\end{eqnarray*}
Due to \eqref{6}, it holds that
\begin{eqnarray*}
\partial_j\eta^i(x,t)=\delta_j^i+\int_0^t\partial_jv^i(x,s)d s.
\end{eqnarray*}
Thus, $D\eta$ can be regarded as a small perturbation of the identity matrix, which implies both $D\eta$ and $A$ are positive definite matrices. Thereby, there exist two positive numbers $\Lambda_1\leq\Lambda_2$ such that
\begin{eqnarray}\label{53}
\Lambda_1|\xi|^2\leq b_k^ib_k^j\xi_j\xi_i\leq\Lambda_2|\xi|^2\quad\ for\ all\ \xi\in \mathbb{R}^n\quad\ and\quad (x,t)\in \Omega\times (0,T^*].
\end{eqnarray}
It follows from the definition of cofactor matrices that
\begin{eqnarray*}
{|B_i^j|\leq (1+M T)^{n-1}}.
\end{eqnarray*}
Note that (see \cite{MB})
\begin{eqnarray*}
J_t=J\textrm{div}u.
\end{eqnarray*}
The chain rule gives
\begin{eqnarray*}
J_t=JB_i^j\partial_jv^i=b_i^j\partial_jv^i.
\end{eqnarray*}
Taking a positive time $T<T^*$  sufficiently small such that { $T\leq\frac{1}{2^{n+1}M}$}, then one has
{ \begin{eqnarray*}
 |J(x,t)-1|&=&|\int_0^tb_i^j(x,s)\partial_jv^i(x,s)d s|\leq J T|B_i^j||\partial_jv^i|\nonumber\\
&\leq& J M T(1+M T)^{n-1}\leq\frac{J}{4},\ in\ \Omega\times(0,T].
\end{eqnarray*}}
This implies
\begin{eqnarray}\label{54}
\frac{1}{2}< J(x,t)<\frac{3}{2}, \ in\ \Omega\times(0,T].
\end{eqnarray}
Direct calculations show (see also \cite{CLS})
\begin{eqnarray*}
\partial_iJ=b_k^j\partial_{ij}\eta^k
\end{eqnarray*}
and
\begin{eqnarray*}
\partial_j{b_i^k}=J^{-1}\partial_{sj}\eta^r(b_r^sb_i^k-b_i^sb_r^k).
\end{eqnarray*}
Therefore, one gets that
{ \begin{eqnarray}\label{55}
|\partial_iJ|\leq \frac{3}{2}(1+M T)^{n-1}MT
\end{eqnarray}}
and
{ \begin{eqnarray}\label{56}
|\partial_j{b_i^k}|\leq 9 (1+M T)^{2n-2}MT.
\end{eqnarray}}

Thus, the system \eqref{52}  is a well-defined integro-differential system with a degeneracy  for t-derivative since the initial density $\rho_0$ vanishes on the boundary $\partial \Omega$.

Define the  linear parabolic operator  $\rho_0\partial_t+L$ by
\begin{eqnarray*}
\rho_0\partial_tw+Lw:&=&\rho_0\partial_tw-\frac{\kappa(\gamma-1) }{R}J^{-1}b_k^ib_k^j\partial_{ij}w\nonumber\\
&&-\frac{\kappa(\gamma-1) }{R}b_k^i\partial_i(J^{-1}b_k^j)\partial_jw
+(\gamma-1)J^{-1}\rho_0 b_i^j\partial_jv^iw.
\end{eqnarray*}
Then,  it  follows from the second equation of \eqref{52} that
\begin{eqnarray}\label{57}
\rho_0\partial_t\mathfrak{e}+L\mathfrak{e}=\frac{\mu}{2}J^{-1}|b_l^j\partial_jv^i+(b_l^j\partial_jv^i)^*|^2+\lambda J^{-1}(b_i^j\partial_jv^i)^2.
\end{eqnarray}

In the rest of this section, our main task is to establish the Hopf's lemma and a strong maximum principle for solutions of  the following  differential inequality
\begin{eqnarray}\label{58}
\rho_0\partial_tw+Lw
\geq0, \ in \ \Omega\times (0,T].
\end{eqnarray}
It follows from \eqref{57} and \eqref{49.4} that $\mathfrak{e}$  also satisfies \eqref{58}.

 We first derive a weak maximum principle for the
  differential inequality \eqref{58}.
\begin{Lemma}\label{1019}
Suppose that $w\in C^{2,1}(Q_T)\cap C(\bar{Q}_T)$ satisfies \eqref{58}.
If $w\geq0(>0)$ on $\partial_p Q_T$,
then $w\geq0(>0)$ in $Q_T$.
\end{Lemma}
\textbf{Proof.} Set %{\red (ÏÂÃæÊÇÐÁÀÏʦдµÄÁíÒ»ÖÖÐÎʽ£¬ÓëÎÒÃÇ×î³õµÄÐÎʽÊǵȼ۵Ä(°´¶¨Òå¿É¿´³ö)£¬ÕâÀï²ÉÓÃÁËËûµÄд·¨£¬¸ü¼òµ¥£¬µ«²»ÈÝÒ׿´³öºóÃæµÄÔËËã)}
\begin{eqnarray*}
d={(\gamma-1)\max_{\bar{\Omega}\times[0,T]}|\frac{J_t}{J}|}
\end{eqnarray*}
and
\begin{eqnarray*}
\varphi=\exp(dt)w.
\end{eqnarray*}
Define a new linear parabolic operator by
\begin{eqnarray*}
\rho_0\partial_t\varphi+\tilde{L}\varphi:=\rho_0\partial_t\varphi+\tilde{L}\varphi
-d\rho_0\varphi.
\end{eqnarray*}
Direct calculation shows that
\begin{eqnarray*}
\rho_0\partial_t\varphi+\tilde{L}\varphi=\exp(dt)(\rho_0\partial_tw+Lw)\geq0,\ in\ Q_T.
\end{eqnarray*}

We first prove the statement under a stronger hypothesis than \eqref{58} that
\begin{eqnarray}\label{59}
\rho_0\partial_t\varphi+\tilde{L}\varphi>0,\ in\ Q_T.
\end{eqnarray}
Assume that $\varphi$ attains its  non-negative minimum at an interior point $(x_0,t_0)$ of the domain $Q_T$. Therefore
\begin{eqnarray*}
\partial_t\varphi(x_0,t_0)\leq0,  \partial_j\varphi(x_0,t_0)=0, a_k^ia_k^j\partial_{ij}\varphi(x_0,t_0)\geq0,
\end{eqnarray*}
which implies $\rho_0\partial_t\varphi+L\varphi\leq0$, this contradicts   \eqref{59}. Next,  choose the auxiliary function
\begin{eqnarray*}
\psi^\varepsilon=\varphi+\varepsilon t,
\end{eqnarray*}
for a  positive number $\varepsilon$. One calculates
\begin{eqnarray*}
\rho_0\partial_t\psi^\varepsilon+L\psi^\varepsilon&=&\rho_0\varphi_t+L\varphi+\varepsilon \rho_0>0,\ in\ Q_T.
\end{eqnarray*}
Thus $\psi^\varepsilon$ attains its  non-negative minimum on $\partial_p Q_T$, which implies that $\varphi$ also attains its  non-negative minimum on $\partial_p Q_T$ by letting $\varepsilon$ go to zero.

Since $w\geq0(>0)$ on $\partial_p Q_T$, so
$\varphi\geq0(>0)$ on $\partial_p Q_T$ by the definition of $\varphi$, furthermore, $\varphi\geq0(>0)$ on $Q_T$. Therefore, $w\geq0(>0)$ on $Q_T$.\qed

The result in Lemma \ref{1019} can also be extended to a general  domain $D\subset \Omega\times (0,T]$.
\begin{Lemma}\label{1020}
Suppose that $w\in C^{2,1}(D)\cap C(\bar{D})$ satisfies \eqref{58}.
If $w\geq0(>0)$ on $\partial_p D$,
then $w\geq0(>0)$ in $D$.
\end{Lemma}

Next, we establish the  Hopf's  lemma for the
  differential inequality \eqref{58}, which is critical for proving  Theorem  \ref{1018}.
\begin{Proposition}\label{1021}
Suppose that $w\in C^{2,1}(\Omega\times (0,T])\cap C(\bar{\Omega}\times [0,T])$ satisfies \eqref{58}
and there exits a point $(x_0,t_0)\in \partial\Omega\times (0,T]$ such that $w(x,t)>w(x_0,t_0)$ for any point $(x,t)$ in  $D$, where %{\red (ÐÁÀÏʦָ³ö ÐèÒª$\Omega$ ±ß½çµÄ¹â»¬ÐÔ£¬ÒÑÔÚÒýÑÔ±ê³ö)}
\begin{eqnarray*}
D=:\{(x,t): |x-\tilde{x}|^2+(t_0-t)< r^2,0< |x-x_0|<\frac{r}{2},0<t\leq t_0\}
\end{eqnarray*}
with $|x_0-\tilde{x}|=r$ and $(x_0-\tilde{x})\perp\partial\Omega$ at $x_0$.
 Then it holds that
\begin{eqnarray*}
\frac{\partial w(x_0,t_0)}{\partial \vec{n}}<0,
\end{eqnarray*}
where $\vec{n}=\frac{x_0-\tilde{x}}{|x_0-\tilde{x}|}$.
\end{Proposition}
\textbf{Proof.}
For positive constants $\alpha$ and $\varepsilon$ to be determined,  set
\begin{eqnarray*}
q(\alpha,x,t)=-e^{-\alpha[|x-\tilde{x}|^2+(t_0-t)]}+e^{-\alpha r^2}
\end{eqnarray*}
and
\begin{eqnarray*}
\varphi(\varepsilon,\alpha,x,t)=w(x,t)-w(x_0,t_0)+\varepsilon q(\alpha,x,t).
\end{eqnarray*}

First, we determine  $\varepsilon$.
The parabolic boundary $\partial_p D$ consists of two parts $\Sigma_1$ and $\Sigma_2$ given by
\begin{gather*}
\Sigma_1=\{(x,t):|x-\tilde{x}|^2+(t_0-t)<r^2, |x-x_0|=\frac{r}{2}, 0<t\leq t_0\}
\end{gather*}
and
\begin{gather*}
\Sigma_2=\{(x,t):|x-\tilde{x}|^2+(t_0-t)=r^2, 0\leq |x-x_0|\leq\frac{r}{2}, 0<t\leq t_0\}.
\end{gather*}
{On  $\bar{\Sigma}_1$, $w(x,t)-w(x_0,t_0)>0$ }, and hence $w(x,t)-w(x_0,t_0)>\varepsilon_0$ for some $\varepsilon_0>0$.
Note that $q\geq-1$ on $\Sigma_1$. Then for such an $\varepsilon_0$,  $\varphi(\varepsilon_0,\alpha,x,t)>0$ on $\Sigma_1$. For $(x,t)\in \Sigma_2$, $q=0$ and $w(x,t)-w(x_0,t_0)\geq0$. Thus, $\varphi(\varepsilon_0,\alpha,x,t)\geq0$
for any $(x,t)\in \Sigma_2$ and $\varphi(\varepsilon_0,\alpha,x_0,t_0)=0$. One concludes that
\begin{eqnarray}\label{60}
\left\{ \begin{array}{ll}
\varphi(\varepsilon_0,\alpha,x,t)\geq0,\  on\  \partial_p D,\\
\varphi(\varepsilon_0,\alpha,x_0,t_0)=0.
\end{array}
\right.
\end{eqnarray}

Next, we choose  $\alpha$. In view of \eqref{58}, one has
\begin{eqnarray}\label{61}
&&\rho_0\partial_t\varphi(\varepsilon_0,\alpha,x,t)+L\varphi(\varepsilon_0,\alpha,x,t)\nonumber\\
&=&\rho_0\partial_tw(x,t)+Lw(x,t)
+\varepsilon_0[\rho_0\partial_tq(\alpha,x,t)+L q(\alpha,x,t)]\nonumber\\
&\geq&\varepsilon_0[\rho_0\partial_tq(\alpha,x,t)+L q(\alpha,x,t)].
\end{eqnarray}
A direct calculation yields
\begin{eqnarray}\label{62}
&&e^{\alpha[|x-\tilde{x}|^2+(t_0-t)]}[\rho_0\partial_tq(\alpha,x,t)+L q(\alpha,x,t)]\nonumber\\
&=&\frac{4\kappa(\gamma-1) }{R}J^{-1}b_k^ib_k^j(x_i-\tilde{x}_i)(x_j-\tilde{x}_j)\alpha^2
-[\rho_0+\frac{2\kappa(\gamma-1) }{R}J^{-1}b_k^ib_k^j\delta_{ij}\nonumber\\
&&+\frac{2\kappa(\gamma-1) }{R}b_k^i\partial_i(J^{-1}b_k^j)(x_j-\tilde{x}_j)]\alpha
-(\gamma-1)J^{-1}\rho_0 b_i^j\partial_jv^i\nonumber\\
&&\times(1-e^{\alpha[|x-\tilde{x}|^2+(t_0-t)- r^2]}).
\end{eqnarray}
It follows from \eqref{53} and \eqref{54} that
\begin{eqnarray}\label{63}
&&\frac{4\kappa(\gamma-1) }{R}J^{-1}b_k^ib_k^j(x_i-\tilde{x}_i)(x_j-\tilde{x}_j)\nonumber\\
&\geq&\frac{8\kappa(\gamma-1)\Lambda_1 }{R}(|x_0-\tilde{x}|-|x-x_0|)^2
\geq\frac{2\kappa(\gamma-1)r^2\Lambda_1}{R}.
\end{eqnarray}
The other terms on the right hand side of \eqref{62} can be estimated by \eqref{55} and \eqref{56} as follows
{ %\red
\begin{align}\label{64}
|\frac{2\kappa(\gamma-1) }{R}J^{-1}b_k^ib_k^j\delta_{ij}|\leq& \frac{4\kappa(\gamma-1)\Lambda_2 }{R},
\\
%\begin{aligned}
|\frac{2\kappa(\gamma-1) }{R}b_k^i\partial_i(J^{-1}b_k^j)(x_j-\tilde{x}_j)|
\leq& \frac{81\kappa(\gamma-1)r }{R}(1+M T)^{3n-3}M T\nonumber\\
\leq& \frac{81\cdot2^{2n-4}\kappa(\gamma-1)r }{R},  \label{65}
%\end{aligned}
\\
% \begin{aligned}
|(\gamma-1)J^{-1}\rho_0 b_i^j\partial_jv^i(1-e^{\alpha[|x-\tilde{x}|^2+(t_0-t)- r^2]})|
\leq& 3(\gamma-1)M^2 (1+M T)^{n-1}\nonumber\\
\leq& 3\cdot2^{n-1}(\gamma-1)M^2,\label{66}
%\end{aligned}
\end{align}}
where \eqref{53}-\eqref{56} have been used.
Finally, one gets
{ %\red
\begin{eqnarray*}
&&e^{\alpha[|x-\tilde{x}|^2+(t_0-t)]}[\rho_0\partial_tq(\alpha,x,t)+L q(\alpha,x,t)]\nonumber\\
&\geq& \frac{2\kappa(\gamma-1)r^2\Lambda_1 }{R}\alpha^2-(M+\frac{4\kappa(\gamma-1)\Lambda_2 }{R}+\frac{81\cdot2^{2n-4}\kappa(\gamma-1)r }{R})\alpha- 3\cdot2^{n-1}(\gamma-1)M^2.
\end{eqnarray*}}
Thereby, there exists a positive number $\alpha_0=\alpha_0(\kappa,\gamma,r,R,M,\Lambda_1,\Lambda_2)$ such that
\begin{eqnarray}\label{66.5}
\rho_0\partial_tq(\alpha_0,x,t)+L q(\alpha_0,x,t)\geq0,\ in \ D.
\end{eqnarray}

 In conclusion, in view of  \eqref{60}, \eqref{61} and \eqref{66.5}, one has
\begin{eqnarray}\label{67}
\left\{ \begin{array}{ll}
\rho_0\partial_t\varphi(\varepsilon_0,\alpha_0,x,t)+L\varphi(\varepsilon_0,\alpha_0,x,t)\geq 0,\ in \ D,\\
\varphi(\varepsilon_0,\alpha_0,x,t)\geq0, \ on \ \partial_pD,\\
\varphi(\varepsilon_0,\alpha_0,x_0,t_0)=0.
\end{array}
\right.
\end{eqnarray}
Lemma \ref{1020}, together with  \eqref{67}, shows that
\begin{eqnarray*}
\varphi(\varepsilon_0,\alpha_0,x,t)\geq0, \ in\ D.
\end{eqnarray*}
Therefore,  $\varphi(\varepsilon_0,\alpha_0,\cdot,\cdot)$ attains its minimum at the point $(x_0,t_0)$ in $D$. In particular, it holds that
\begin{eqnarray*}
\varphi(\varepsilon_0,\alpha_0,x,t_0)\geq \varphi(\varepsilon_0,\alpha_0,x_0,t_0)\quad  \ for\ all\  x\in \{x:|x-x_0|\leq\frac{r}{2}\}.
\end{eqnarray*}
This  implies
\begin{eqnarray*}
\frac{\partial \varphi(\varepsilon_0,\alpha_0,x_0,t_0)}{\partial \vec{n}}\leq0.
\end{eqnarray*}
Finally, one obtains
\begin{eqnarray*}
\frac{\partial w(x_0,t_0)}{\partial \vec{n}}\leq-\varepsilon_0 \frac{\partial q(\alpha_0,x_0,t_0)}{\partial \vec{n}}
=-2\varepsilon_0\alpha_0 r e^{-\alpha_0 r^2}<0.
\end{eqnarray*}
\qed

In order to establish a strong maximum principle for the
  differential inequality \eqref{58}, we study first the t-derivative at an
interior minimum point.
\begin{Lemma}\label{1022}
Let  $w\in C^{2,1}(\Omega\times (0,T])\cap C(\bar{\Omega}\times [0,T])$ satisfy \eqref{58} and have a minimum $M_0$ in the domain $\Omega\times (0,T]$. Suppose that $\Omega\times (0,T]$ contains a closed solid ellipsoid
\begin{eqnarray*}
\Omega^\sigma:=\{(x,t): |x-x_*|^2+\sigma(t-t_*)^2\leq r^2\},\sigma>0
\end{eqnarray*}
and  $w(x,t)>M_0$ for any interior point $(x,t)$ of $\Omega^\sigma$ and  $w(\bar{x},\bar{t})=M_0$ at some point $(\bar{x},\bar{t})$ on the  boundary   of $\Omega^\sigma$. Then $\bar{x}=x_*$.
\end{Lemma}
\textbf{Proof.}  One can assume that $(\bar{x},\bar{t})$ is the only point on $\partial \Omega^\sigma$ such that $w=M_0$ in $\Omega^\sigma$.  Otherwise,  one can limit it to a smaller closed ellipsoid in $\Omega^\sigma$ and with $(\bar{x},\bar{t})$ as the only common point with $\partial \Omega^\sigma$. We prove the desired result by contradiction.
Suppose that $\bar{x}\neq x_*$.  Choose a  closed ball $D$ with center $(\bar{x},\bar{t})$ and radius $\tilde{r}<|\bar{x}-x_*|$ contained in $\Omega\times (0,T]$. Then, one has
\begin{eqnarray}\label{68}
|x-x_*|\geq|\bar{x}-x_*|-\tilde{r}=:\hat{r}\quad \ for\  (x,t)\in D.
\end{eqnarray}
The parabolic boundary $\partial_pD=\partial D$ of $D$ consists of a part $\Sigma_1$ lying in $\Omega^\sigma$ and a part $\Sigma_2$ lying outside $\Omega^\sigma$.

For positive constants $\alpha$ and $\varepsilon$ to be determined, set
\begin{eqnarray*}
q(\alpha,x,t)=-e^{-\alpha[|x-x_*|^2+\sigma(t-t_*)^2]}+e^{-\alpha r^2}
\end{eqnarray*}
and
\begin{eqnarray*}
\varphi(\varepsilon,\alpha,x,t)=w(x,t)-M_0+\varepsilon q(\alpha,x,t).
\end{eqnarray*}
We first determine the value of $\varepsilon$. Note that $q(\alpha,x,t)<0$ in the interior of $\Omega^\sigma$, $q(\alpha,x,t)=0$ on $\partial \Omega^\sigma$ and $q(\alpha,x,t)>0$ outside $\Omega^\sigma$. So, it holds that $\varphi(\varepsilon,\alpha,\bar{x},\bar{t})=0$.
On  $\Sigma_1$, $w(x,t)-M_0>0$, and hence $w(x,t)-M_0>\varepsilon_0$ for some $\varepsilon_0>0$.
Note that $q(\alpha,x,t)\geq-1$ on $\Sigma_1$. Then for such an $\varepsilon_0$, $\varphi(\varepsilon_0,\alpha,x,t)>0$ on $\Sigma_1$. {For $(x,t)\in \Sigma_2$, we have $q(\alpha,x,t)>0$ and $w(x,t)-M_0\geq0$}. Thus, $\varphi(\varepsilon_0,\alpha,x,t)>0$
for any $(x,t)\in \Sigma_2$. One concludes that
\begin{eqnarray}\label{69}
\left\{ \begin{array}{ll}
\varphi(\varepsilon_0,\alpha,x,t)>0,\  on\  \partial_pD,\\
\varphi(\varepsilon_0,\alpha,\bar{x},\bar{t})=0.
\end{array}
\right.
\end{eqnarray}

Next, we choose  $\alpha$.
We  need to estimate $\rho_0q_t(\alpha,x,t)+L q(\alpha,x,t)$ due to \eqref{61}.
One calculates
\begin{eqnarray*}
&&e^{\alpha[|x-x_*|^2+\sigma(t-t_*)^2]}[\rho_0\partial_tq(\alpha,x,t)+L q(\alpha,x,t)]\nonumber\\
&=&\frac{4\kappa(\gamma-1) }{R}J^{-1}b_k^ib_k^j(x_i-(x_*)_i)(x_j-(x_*)_j)\alpha^2
-[2\sigma\rho_0(t-t_*)\nonumber\\
&&+\frac{2\kappa(\gamma-1) }{R}J^{-1}b_k^ib_k^j\delta_{ij}+\frac{2\kappa(\gamma-1) }{R}b_k^i\partial_i(J^{-1}b_k^j)(x_j-{x_*}_j)]\alpha\nonumber\\
&&-(\gamma-1)J^{-1}\rho_0 b_i^j\partial_jv^i(1-e^{\alpha[|x-x_*|^2+\sigma(t-t_*)^2- r^2]}).
\end{eqnarray*}
Similar to \eqref{63}-\eqref{66},  there exists a positive number $\alpha_0=\alpha_0(\kappa,\gamma,\sigma,r,\hat{r},R,M,\Lambda_1,\Lambda_2)$ such that
\begin{eqnarray}\label{70}
\rho_0\partial_tq(\alpha_0,x,t)+L q(\alpha_0,x,t)\geq0,\ in \ D.
\end{eqnarray}

 In conclusion, it follows  from \eqref{61} and \eqref{70} that
\begin{eqnarray}\label{71}
\left\{ \begin{array}{ll}
\rho_0\partial_t\varphi(\varepsilon_0,\alpha_0,x,t)+L\varphi(\varepsilon_0,\alpha_0,x,t)\geq 0,\ in \ D,\\
\varphi(\varepsilon_0,\alpha_0,x,t)>0, \ on \ \partial_pD,\\
\varphi(\varepsilon_0,\alpha_0,\bar{x},\bar{t})=0.
\end{array}
\right.
\end{eqnarray}
Then Lemma \ref{1020} and \eqref{71} imply that
\begin{eqnarray*}
\varphi(\varepsilon_0,\alpha_0,x,t)>0, \ in\ D.
\end{eqnarray*}
which contradicts $\varphi(\varepsilon_0,\alpha_0,\bar{x},\bar{t})=0$ due to $(\bar{x},\bar{t})\in D$ .\qed

Based on Lemma \ref{1022}, it is standard to prove the following lemma.  For details,  one can refer to Lemma 3 of Chapter 2 in \cite{Friedman}.
\begin{Lemma}\label{1023}
Suppose that $w\in C^{2,1}(\Omega\times (0,T])\cap C(\bar{\Omega}\times [0,T])$ satisfies \eqref{58}.
If $w$ has a minimum in an interior point $P_0=(x_0,t_0)$  of  $\Omega\times (0,T]$, then
$w(P)=w(P_0)$ for any point $P(x,t_0)$ of $\Omega\times (0,T]$.
\end{Lemma}

Next, we prove a local  strong minimum principle in a rectangle $\mathcal{R}$  of the domain $\Omega\times (0,T]$.
\begin{Lemma}\label{1024}
Suppose that $w\in C^{2,1}(\Omega\times (0,T])\cap C(\bar{\Omega}\times [0,T])$ satisfies \eqref{58}.
If $w$ has a minimum in the interior point $P_0=(x_0,t_0)$ of $\Omega\times (0,T]$, then there exists a rectangle
\begin{eqnarray*}
\mathcal{R}(P_0):=\{(x,t):(x_0)_i-c_i\leq x_i\leq (x_0)_i+c_i, t_0-c_0\leq t\leq t_0,i=1,2,\cdot\cdot\cdot,n\}
\end{eqnarray*}
in $\Omega\times (0,T]$ such that
$w(P)=w(P_0)$ for any point $P$ of $\mathcal{R}(P_0)$.
\end{Lemma}
\textrm{\textbf{Proof.}} We prove the desired result by contradiction. Suppose that there exists an interior point $P_1=(x_1,t_1)$ of $\Omega\times (0,T]$ with $t_1< t_0$ such that $w(P_1)>w(P_0)$.
 Connect  $P_1$ to $P_0$ by a simple smooth curve $\gamma$. Then there exists a point $P_*=(x_*,t_*)$ on $\gamma$ such that $w(P_*)=w(P_0)$ and $w(\bar{P})<w(P_*)$ for all any point $\bar{P}$ of $\gamma$ between $P_1$ and $P_*$. We may assume that $P_*=P_0$ and $P_1$ is very near to $P_0$. There exists a rectangle $\mathcal{R}(P_0)$ in $\Omega\times (0,T]$ with small positive numbers $a_0$ and $a_1$ (to be determined) such that
 $P_1$ lies on $t=t_0-a_0$.
  Since $\mathcal{R}(P_0)\setminus \{t=t_0\}\cap\{t=\bar{t}\}$ contains some point $\bar{P}=(\bar{x},\bar{t})$ of $\gamma$ and $w(\bar{P})>w(P_0)$,  we deduce $w(P)>w(P_0)$ for each point $P$ in $\mathcal{R}(P_0)\setminus \{t=t_0\}\cap\{t=\bar{t}\}$ due to Lemma  \ref{1011}.  Therefore, $w(P)>w(P_0)$ for each point $P$ in $\mathcal{R}(P_0)\setminus \{t=t_0\}$.

For positive constants $\alpha$ and $\varepsilon$ to be determined, set
\begin{eqnarray*}
q(\alpha,x,t)=-t_0+t+\alpha|x-x_0|^2
\end{eqnarray*}
and
\begin{eqnarray*}
\varphi(\varepsilon,\alpha,x,t)=w(x,t)-w(P_0)+\varepsilon q(\alpha,x,t).
\end{eqnarray*}
Assume further that $P=(x_0-c,t_0-c_0)$ is on the parabola $q(\alpha,x,t)=0$, then
\begin{equation}\label{71.2}
  \alpha=\frac{c_0}{|c|^2},
\end{equation}
where $|c|=(\sum_{i=1}^n|c_i|^2)^{\frac{1}{2}}$.

A direct calculation shows that
\begin{eqnarray}\label{71.4}
&&\rho_0\partial_tq(\alpha,x,t)+L q(\alpha,x,t)\nonumber\\
&=&-\alpha[\frac{2\kappa(\gamma-1) }{R}J^{-1}b_k^ib_k^j\delta_{ij}
+\frac{2\kappa(\gamma-1) }{R}b_k^i\partial_i(J^{-1}b_k^j)(x_j-(x_0)_j)\nonumber\\
&&-(\gamma-1)J^{-1}\rho_0 b_i^j\partial_jv^i|x-x_0|^2)]+\rho_0[1+(\gamma-1)J^{-1} b_i^j\partial_jv^i(-t_0+t)].
\end{eqnarray}
The first three terms on the right hand side of \eqref{73} can be estimated similar to \eqref{63}-\eqref{65}. For the last term, one has
{ %\red
\begin{gather*}
|(\gamma-1)J^{-1} b_i^j\partial_jv^i(-t_0+t)|
\leq 6(\gamma-1)(1+M T)^{n-1}MT\leq \frac{3}{2}(\gamma-1).
\end{gather*}}
Consequently, one gets
{ %\red
\begin{eqnarray}\label{72}
&&\rho_0\partial_tq(\alpha,x,t)+L q(\alpha,x,t)\nonumber\\
&\geq&-\alpha[\frac{\kappa(\gamma-1) }{R}(4\Lambda_2+81\cdot2^{2n-4}|c|)+3\cdot2^{n-1}(\gamma-1)M^2|c|^2]+\frac{3\gamma-1}{2}\rho_0.
\end{eqnarray}}
Since $\rho_0$ has a positive lower bound depending on $x_0\pm c$ in $\mathcal{R}(P_0)$, one can choose $\alpha_0$ such that
{ %\red
\begin{eqnarray}\label{73}
\alpha_0 <\frac{(3\gamma-1)R\rho_0}{\kappa(\gamma-1)(8\Lambda_2+81\cdot2^{2n-3}|c|) +3\cdot2^n(\gamma-1)RM^2|c|^2},
\end{eqnarray}}
then it follows from \eqref{71.4}-\eqref{73} that
 \begin{eqnarray}\label{74}
\rho_0\partial_t\varphi(\alpha_0,x,t)+L \varphi(\alpha_0,x,t)
\geq0, \ in\ \mathcal{R}(P_0).
\end{eqnarray}

Next, for the fixed $c_0$, one can choose $c$ such that $\mathcal{R}(P_0)\subset\Omega\times (0,T]$ and then it follows from \eqref{71.2} and \eqref{74} that $c_0$ can be choosen such that
{ %\red
 \begin{equation*}
  c_0<\min\{t_0,\frac{(3\gamma-1)|c|^2R\rho_0}{\kappa(\gamma-1)(16\Lambda_2+81\cdot2^{2n-2}|c|) +3\cdot2^{n+1}(\gamma-1)RM^2|c|^2}\}.
\end{equation*}}

Denote $\mathcal{S}=\{(x,t)\in \mathcal{R}(P_0), q(x,t)\geq0\}$.  The parabolic boundary $\partial_p\mathcal{S}$ of $\mathcal{S}$ consists of a part $\Sigma_1$ lying in $\mathcal{R}(P_0)$ and a part $\Sigma_2$ lying on  $\mathcal{R}(P_0)\cap \{t=t_0-c_0\}$.

Finally, one can choose  $\varepsilon$.
On  $\Sigma_2$, $w(x,t)-M_0>0$. Note $q(\alpha,x,t)$ is bounded on $\Sigma_2$, one can choose $\varepsilon_0$ suitably small such that $\varphi(\varepsilon_0,\alpha_0,x,t)>0$ on $\Sigma_2$.
 On $\Sigma_1\setminus \{P_0\}$, $q(\alpha,x,t)=0$ and $w(x,t)-M_0>0$. Thus, $\varphi(\varepsilon_0,\alpha_0,x,t)>0$ on $\Sigma_1\setminus \{P_0\}$ and
 $\varphi(\varepsilon_0,\alpha_0,x_0,t_0)=0$. One concludes that
\begin{eqnarray}\label{75}
\left\{ \begin{array}{ll}
\varphi(\varepsilon_0,\alpha_0,x,t)>0,\  on\  \partial_p\mathcal{S}\setminus \{P_0\},\\
 \varphi(\varepsilon_0,\alpha_0,x_0,t_0)=0.
\end{array}
\right.
\end{eqnarray}

 In conclusion,  it follows from \eqref{74} and \eqref{75} that
\begin{eqnarray}\label{76}
\left\{ \begin{array}{ll}
\rho_0\partial_t\varphi(\varepsilon_0,\alpha_0,x,t)+L\varphi(\varepsilon_0,\alpha_0,x,t)\geq 0,\ in\ \mathcal{S},\\
\varphi(\varepsilon_0,\alpha_0,x,t)>0,\  on\  \partial_p\mathcal{S}\setminus \{P_0\},\\
 \varphi(\varepsilon_0,\alpha_0,x_0,t_0)=0.
\end{array}
\right.
\end{eqnarray}
In view of Lemma \ref{1020} and \eqref{76},  $\varphi(\varepsilon_0,\alpha_0,\cdot,\cdot)$  attains its minimum at $P_0$ in $\mathcal{\bar{S}}$, thus
\begin{eqnarray*}
\frac{\partial \varphi(\varepsilon_0,\alpha_0,x_0,t_0)}{\partial t}\leq0.
\end{eqnarray*}
Note that $q$ satisfies at $P_0$
\begin{eqnarray*}
\frac{\partial q(\alpha_0,x_0,t_0)}{\partial t}=1.
\end{eqnarray*}
 Therefore
\begin{eqnarray}\label{77}
\frac{\partial w(x_0,t_0)}{\partial t}\leq-\varepsilon_0.
\end{eqnarray}
But, by the assumption, $w$  attains its   minimum   at $P_0$,  it follows that
\begin{eqnarray*}
\rho_0\frac{\partial w(x_0,t_0)}{\partial t}\geq -L w(x_0,t_0)\geq0,
\end{eqnarray*}
which contradicts to \eqref{77}.\qed

 Now the following global strong maximum principle can be proved similarly as for Proposition \ref{1013}.
\begin{Proposition}\label{1025}
Suppose that $w\in C^{2,1}(\Omega\times (0,T])\cap C(\bar{\Omega}\times [0,T])$ satisfies \eqref{58}.
If $w$ attains its minimum at some interior point $P_0=(x_0,t_0)$ of $ \Omega\times (0,T]$, then
$w(P)=w(P_0)$ for any point  $P$ of $\Omega\times (0,t_0]$.
\end{Proposition}

We are ready to prove Theorem \ref{1018}.\\
\textrm{\textbf{Proof of Theorem \ref{1018}.}}
Recall that $\mathfrak{e}$ satisfies \eqref{58}, so the weak maximum principle, Hopf lemma  and strong maximum principle holds for $\mathfrak{e}$.
Since  $\mathfrak{e}_0\geq0$ and $\mathfrak{e}_0\not\equiv0$ in $ \Omega$, and $\mathfrak{e}=0$ on  $\partial\Omega\times (0,t_0]$ due to \eqref{52}, by Proposition \ref{1025}, it holds that $\mathfrak{e}>0$ in $\Omega\times(0,T]$.   Taking any point $(x_0,t_0)$ of $\partial\Omega\times(0,T]$, applying Proposition \ref{1021}, we obtain $\frac{\partial \mathfrak{e}(x_0,t_0)}{\partial \vec{n}}<0$, which contradicts to $\mathfrak{e}_{x_i}(x_0,t_0)=0$ on $\partial\Omega\times (0,T]$ due to \eqref{52}. \qed

\vskip 5mm

\noindent \textrm{\textbf{Acknowledgements:}} The research of Li was supported partially by the National Natural Science Foundation of China (Nos. 11231006, 11225102, 11461161007 and 11671384), and the Importation and Development of High Caliber Talents Project of Beijing Municipal Institutions (No. CIT\&TCD20140323).
The research of Wang was supported by grant nos. 231668 and 250070 from the Research Council of Norway.
 The research of Xin was supported partially by the Zheng Ge Ru Foundation, Hong Kong RGC Earmarked Research grants CUHK-14305315 and CUHK-4048/13P, NSFC/RGC Joint Research Scheme N-CUHK443/14, and Focused Innovations Scheme from The Chinese University of Hong Kong.

%\end{CJK*}

\textsc{School of Mathematics, Capital Normal University, Beijing 100048, P. R. China}

 E-mail address: hailiang.li.math@gmail.com

\textsc{Department of Mathematical Sciences, Norwegian University of Science and Technology, Trondheim 7491, Norway}

E-mail address: yuexun.wang@ntnu.no

\textsc{The Institute of Mathematical Sciences, The Chinese University of Hong Kong, Hong Kong}

E-mail address: zpxin@ims.cuhk.edu.hk

\end{document}